\documentclass{siamltex}

\usepackage{amsmath}
\usepackage{amssymb}

\usepackage[toc,page]{appendix}

\usepackage{tikz}
\usepackage{pgfplots}
\usepgfplotslibrary{groupplots}
\usepackage[ComplainAboutUnexternalized]{externaltikz}
\usepackage{calc}

\usepackage{subcaption}
\usepackage{hyperref}
 \usepackage{mathtools}

\newtheorem{remark}{Remark}
\newtheorem{example}{Example}

\newcommand{\R}{\mathbb{R}}

\begin{document}
	
\title{Relaxation models for  scalar traffic networks and  zero relaxation limit}

\author{R. Borsche\footnotemark[1] 
	\and  A. Klar\footnotemark[1] \footnotemark[2]}
\footnotetext[1]{Technische Universit\"at Kaiserslautern, Department of Mathematics, Erwin-Schr\"odinger-Stra{\ss}e, 67663 Kaiserslautern, Germany 
	(\{borsche, klar\}@mathematik.uni-kl.de)}
\footnotetext[2]{Fraunhofer ITWM, Fraunhoferplatz 1, 67663 Kaiserslautern, Germany}

	\maketitle

	\begin{abstract}
		In this paper we propose coupling conditions for a relaxation model  for vehicular traffic on networks.  
		We present a matched asymptotic expansion procedure  to derive a LWR- network with    well-known classical coupling conditions
		from the relaxation network in the macroscopic limit. Similar to the asymptotic limit of boundary value problems, we perform 
	 	an asymptotic analysis of the interface  layers at the nodes and  a matching procedure	using half-Riemann problems for the limit conservation law.
	 	Moreover, we present numerical experiments comparing the relaxation network with the LWR network for a broader range of coupling conditions.
	\end{abstract}

\section{Introduction}

Modeling and simulation of traffic flow on road networks, has been investigated intensively using hyperbolic partial differential equations. % to describe the  dynamics on single arcs of the network.
Different models have been used, ranging from scalar conservation laws like the Lighthill Whitham Richards model through 
models using system of conservation laws, \cite{P79, AR,Ber,Hel,G02,Dag1,Dag2}, to kinetic descriptions of the flow \cite{Hel,HPRV20,KW97,PSTV17,FT13}.
Derivations of these models from the underlying models in such hierarchies have been discussed as well, see,  \cite{AKMR,Hel,KW00,B11}, for a non-exhaustive list of references.
To obtain a model for the dynamics on the full network, all these models have to be supplemented with coupling conditions at the nodes of the network. 
Coupling conditions for scalar conservation laws and systems of conservation laws on networks  have been discussed in many papers, see, for example,  \cite{Cor17,BNR14,HR95,GPBook,LS02,CGP05,BHK06a,CG08,CM08,G10,KCGG}.
Kinetic and relaxation equations on networks have been considered, for example,  in  \cite{HM09,BKKP16}. 

A  procedure to derive coupling conditions for macroscopic equations from the underlying  ones of the kinetic or relaxation equation has been  discussed for linear systems in \cite{BK18c} using an asymptotic analysis of the situation near the nodes. A  simple nonlinear case has been treated in \cite{BK18b}.
To explain the general procedure in  more detail, we consider a  relaxation equation in 1D involving a scaling parameter $\epsilon$, which  converges for $\epsilon \rightarrow 0$ to  an associated scalar conservation law for traffic flow.
If such equations are considered on a network, it is sufficient to study a single coupling point or node, where  coupling conditions are required. 
Suitable  coupling conditions have  to be imposed for the relaxation problem at each node.
If $\epsilon$ is send to zero,  layers near  the junctions can arise.
To consider the limit  $\epsilon \rightarrow 0$, one has to proceed similarly
as in the case of  boundary value problems, where  a complete picture of the convergence is only obtained, once  boundary- and initial layers are investigated.
We refer to  \cite{BSS84,BLP79,G08,UTY03} for such a procedure for boundary value problems  in the case of kinetic equations and to  
\cite{WY99,WX99,LX96,X04} for the case of hyperbolic relaxation systems.

In the present work, we consider the case of a relaxation system  on a network with a small parameter $\epsilon$
leading in the limit  $\epsilon \rightarrow 0 $ to a LWR-type scalar conservation law.
Besides the definition of suitable coupling conditions for the relaxation system at the junction,  the present work aims at presenting a matched  asymptotic expansion procedure leading from the  relaxation model on the network to  the  scalar conservation law.
Analytical, as well as numerical investigations are presented.

Based on the discussion of the Riemann problems at the nodes we propose  coupling conditions for the relaxation model for merging and diverging junctions. These conditions are developed in a similiar way as those for other well-known 
higher order traffic   models like the  ARZ-equations, see \cite{HR}. 
However, due to the simpler structure of the relaxation model compared to the ARZ-equations, the conditions are much easier to handle and to investigate and allow a relatively straightforward  asymptotic derivation of classical coupling conditions for the  LWR-type  traffic equations in the limit 
$\epsilon $ going to $0$. Moreover,  in the case of diverging junctions the coupling conditions defined here  guarantee that the coupled solution of the two-equation model remain in the physically reasonable state-space domain defined by bounds on density and velocity in contrast to coupling conditions for the ARZ-model discussed in the literature \cite{HR,KCGG} and references therein.

The asymptotic procedure gives a detailed account of the situation near the node in the case of small values of
the relaxation parameter $\epsilon$.  Besides the structure of the layers near the nodes and the coupling conditions for the limit problem, it reveals, that the effective  densities for the relaxation system at the node for small $\epsilon$  are not necessarily the same as those found by the coupling conditions for the relaxation system.

The paper is organized in the following way.  
In section \ref{equations} we present the  relaxation model, compare \cite{BK18}, and the associated scalar conservation law. 
In section \ref{sec:kineticcouplingconditions} different coupling conditions for the  relaxation model for merging and diverging junctions  are discussed and the associated Riemann problems at the nodes are investigated.
Moreover, the associated classical coupling conditions for the limiting scalar conservation law are stated.
In section \ref{asyproc}  the asymptotic procedure and the matched asymptotic expansion on the network is explained together with a discussion of  the layer solutions of the relaxation system.
Then, the asymptotic procedure  is investigated in detail analytically in section
 \ref{macroscopiccc}. There, it is shown 	that a special  merge conditions for the relaxation system leads in the relaxation  limit to  
 a classical merge condition  for the  nonlinear scalar conservation law.
Finally, the solutions of the relaxation system on the network are compared numerically to the solution of the scalar conservation law  on the network in section \ref{Numerical results} for the case of merging and diverging junctions and a broader range of coupling conditions.

\section{Relaxation model and scalar  traffic equations}
\label{equations}

Consider the LWR-traffic flow equations
\begin{align}\label{LWR0}
\begin{aligned}
\partial_t \rho + \partial_x F(\rho) &=0,\\
\end{aligned}
\end{align}
where $F= F(\rho)$ is a given traffic density-flow function or fundamental diagram, i.e. a smooth function $F:[0,1]\rightarrow[0,1]$ with $F(0)=0=F(1)$ and $F^\prime(\rho) \le 1$.
In the following we restrict ourselves to strictly concave fundamental diagrams $F$, where 
the point, where the maximum of $F$ is attained, is denoted by   $\rho^{\star}$ and the maximal value is 
$ F(\rho^{\star})=\sigma$.

We are interested in the investigation of relaxation systems for the LWR equations on networks. The minimal requirements of such a relaxation system for the density $\rho$ and the flux $q$ are the convergence towards the  LWR-equations and the invariance of the 'traffic domain'
given by  $0 \le \rho \le 1$ and $0 \le q \le \rho$.
This is the physically reasonable state-space for a 2x2 traffic equation, see e.g. \cite{Ber}.
The above two conditions correspond to an upper bound on the density $\rho_{max} =1$ and an upper bound to the 
velocity $v=q/\rho$  given by $v_{max}$.

A simple example is given by the following relaxation system \cite{BK18} for the LWR equations for the variables density $\rho$ and flux $q$, which we will use as a prototype for the discussion of the issues mentioned in the introduction.
The equations are 
\begin{align}\label{macro0}
\begin{aligned}
\partial_t \rho + \partial_x q &=0\\
\partial_t q + \frac{ q}{1-\rho}  \partial_x \rho  + (1-\frac{ q}{1-\rho})\partial_x  q  &=-\frac{1}{\epsilon} \left(q-F(\rho) \right) \ .
\end{aligned}
\end{align}

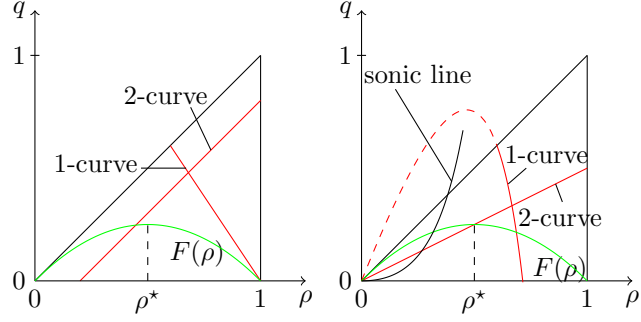
\begin{figure}[h]
	\center
	\externaltikz{Statespacelayer_H1}{
					\begin{tikzpicture}[scale = 3]
		\def\dr{0.2}
		\node[below] at (0,0) {$0$};
		\node[left] at (0,0) {$0$};
		\node[below] at (1,0) {$1$};
		\node[left] at (0,1) {$1$};
		\draw(1,0)--(1,1);
		\draw(0,0)--(1,1);
		\draw(-0.02,1)--(0.02,1);
		
		\node[below] at (0.72,0.22) {$F(\rho)$};
		
		\draw[red](0.2,0.0)--(1,0.8);
		
		\draw(0.78,0.58)--(0.69,0.77);	
		\node[below] at (0.27,0.6) {$1$-curve};	
		
		\draw[red](1.0,0.0)--(0.6,0.6);		
		
		\draw(0.47,0.5)--(0.67,0.5);	
		\node[below] at (0.58,0.9) {$2$-curve};		
		
		\draw[dashed](0.5,0.0)--(0.5,0.25) node[below] at (0.5,0) {$\rho^\star$};	
		\draw[->](0,0)--(1.2,0) node[below]{$\rho$};
		\draw[->](0,0)--(0,1.2) node[left]{$q$};
		\draw[domain=0.0:1,smooth,variable=\x,green] plot ({\x},{\x*(1-\x)});	
		\end{tikzpicture}
		\begin{tikzpicture}[scale = 3]
		\def\dr{0.2}
		\node[below] at (0,0) {$0$};
		\node[left] at (0,0) {$0$};
		\node[below] at (1,0) {$1$};
		\node[left] at (0,1) {$1$};
		\draw(1,0)--(1,1);
		\draw(0,0)--(1,1);
		\draw(-0.02,1)--(0.02,1);
		\node[below] at (0.88,0.15) {$F(\rho)$};
		
		\draw(0.65,0.43)--(0.73,0.5);	
		\node[below] at (0.82,0.64) {$1$-curve};	
		
		\draw[red](0.0,0.0)--(1,0.5);		
		
		\draw(0.86,0.43)--(0.89,0.35);	
		\node[below] at (0.88,0.36) {$2$-curve};		
		
		\draw(0.15,0.85)--(0.4,0.44);	
		\node[below] at (0.25,1.0) {sonic line};	
		\draw[dashed](0.5,0.0)--(0.5,0.25) node[below] at (0.5,0) {$\rho^\star$};	
		\draw[->](0,0)--(1.2,0) node[below]{$\rho$};
		\draw[->](0,0)--(0,1.2) node[left]{$q$};
		\draw[domain=0.0:0.6,smooth,variable=\x,red,dashed] plot ({\x},{\x*(2.5-\x/(1-\x))});
		\draw[domain=0.6:0.715,smooth,variable=\x,red] plot ({\x},{\x*(2.5-\x/(1-\x))});	
		\draw[domain=0.0:1,smooth,variable=\x,green] plot ({\x},{\x*(1-\x)});	
		\draw[domain=0.0:0.45,smooth,variable=\x] plot ({\x},{\x^2/(1-\x)^2});	
		\end{tikzpicture}
	}
	\caption{Lax-curves  in $(\rho,q)$ variables for the  totally degenerate relaxation system (left figure) and the ARZ equations
		with $p(\rho ) = \frac{\rho}{1-\rho}$ (right figure).}
	\label{state0}
\end{figure}

This is a hyperbolic system with
the eigenvalues $\lambda_1 = - \frac{q}{1-\rho}\leq 0 <\lambda_2 = 1$.  
The respective eigenvectors  are
$r_1 = \left(
		1,
		\lambda_1 \right)^T, 
	r_2 = \left(
	1,1
	\right)\ $. 
A straightforward computation shows that the $r_1$- and  the $r_2$-field are both linearly degenerate. 
The system is totally linear degenerate (TLD).
The integral curves (and shock curves) of the hyperbolic system are given by $q= q_L \frac{1-\rho}{1-\rho_L}$ for the 1-field and by $q = \rho -\rho_R+q_R$ for the 2-field. 
The region $0 \le \rho \le 1, 0 \le q \le \rho$ is an invariant region for the  relaxation system as it is easily seen by considering  the 
integral curves, see Figure \ref{state0}.
We refer to \cite{BK18} for details and for a kinetic interpretation of the equations.
Note that the fact, that  the system is TLD strongly simplifies the calculations and leads in many situations to explicitly computable quantities and conditions.

The equations  can be rewritten in conservative form choosing the   variable 
$
	z =  \frac{q}{1-\rho}\ .
$
Rewriting \eqref{macro0} we obtain
	\begin{align}\label{eq:lindeg+relax}
	\begin{aligned}
	\partial_t \rho + \partial_x q&=0\\
	\partial_t z + \partial_x z &= -\frac{1}{\epsilon} \left(z - Z(\rho) \right)\ 
	\end{aligned}
	\end{align}
with $q =  z (1-\rho)$ and $Z(\rho) = \frac{F(\rho)}{1-\rho}$.	
A Riemann invariant of the first characteristic family is obviously 
$$z= \frac{q}{1-\rho} \in [0,\infty)\ .$$
A  Riemann invariant of the  second characteristic family is 
$$w =\rho - q  \in [0,1]\ .$$

Concerning the convergence of its solutions towards the solutions of the scalar conservation law $\partial_t \rho + \partial_x F(\rho) =0$ as
$\epsilon $ tends to $0$ the subcharacteristic condition has to be satisfied \cite{LX96}.
Setting $q = F(\rho)$ in  the formula for the eigenvalues, the subcharacteristic condition states
$$
-  \frac{ F(\rho)}{1-\rho} \le F^\prime(\rho) \le 1 \  \mbox{ for } \  0 \le \rho \le 1\ .
$$

\begin{remark}
	The condition is fulfilled for strictly concave fundamental diagrams $F$.
	For example, in the classical LWR case with $F(\rho) =\rho (1-\rho) $ and  $F^\prime (\rho) = 1- 2 \rho$ the above  condition is
	$$
	-   \rho \le 1- 2 \rho \le 1\  \mbox{ for } \ 0 \le \rho \le 1\ ,
	$$
	which is obviously satisfied. 
\end{remark}

Finally we note that boundary conditions for the relaxation system  \eqref{macro0} on the interval $[x_L,x_R]$ have to be prescribed in the following way. Since the first  eigenvalue is always non-positive and the second is a positive constant,  the number of boundary conditions is fixed. 
At the left boundary at $x=x_L$ we have to prescribe a value for the 1-Riemann invariant $z(x_L) = \frac{q(x_L)}{1-\rho(x_L)} $. For the right boundary $x= x_R$ the 2-Riemann invariant $w(x_R)=\rho(x_R)- q(x_R))$
has to be prescribed.

\begin{remark}
Compare the present model with the ARZ model \cite{AR}	
\begin{align}
\label{macro0rascle}
\partial_t \rho +\partial_x  q &=0\\
\partial_t q  +\frac{q}{\rho}\left(\rho p'(\rho)- \frac{q}{\rho}\right)\partial_x  \rho +\left( \frac{q}{\rho}\ + \frac{q}{\rho}-\rho p'(\rho)  \right)\partial_x q&= -\frac{1}{\epsilon} \left(q-F(\rho) \right). \nonumber
\end{align}
or in conservative form
\begin{align}
\label{macro0konsrascle}
\partial_t \rho +\partial_x  q &=0\\
\partial_t (\rho z_R)   +\partial_x (qz_R)&= -\frac{1}{\epsilon} \left(q-F(\rho) \right). \nonumber
\end{align}
with $q = \rho z_R - \rho p(\rho)$.
We note that the above domain $0 \le \rho \le 1$ and $0 \le q \le \rho$ is also an invariant domain for the ARZ conditions, if $p$ is appropriately chosen
with a singularity at $\rho=1$, see \cite{Ber}.
\end{remark}

\section{Coupling conditions for the relaxation system on networks and associated conditions for the LWR equations} \label{sec:kineticcouplingconditions}
	
\begin{figure}[h!]
	\begin{center}
		\externaltikz{sketch_21node}{
			\begin{tikzpicture}[thick]
			\def\len{2}
			\node[fill,circle] (N) at (0,0){};
			\draw[->] (-\len,0.6)--(N) node[above,pos = 0.5]{$1$};
			\draw[->] (-\len,-0.6)--(N) node[below,pos = 0.5]{$2$};
			\draw[->] (0,0)--(\len,0) node[above,pos = 0.5]{$3$};
			\end{tikzpicture}
		}
	\;\;
		\externaltikz{sketch_12node}{
			\begin{tikzpicture}[thick]
			\def\len{2}
			\node[fill,circle] (N) at (0,0){};
			\draw[->] (-\len,0)--(N) node[above,pos = 0.5]{$1$};
			\draw[->] (0,0)--(\len,0.6) node[above,pos = 0.5]{$2$};
			\draw[->] (0,0)--(\len,-0.6) node[below,pos = 0.5]{$3$};
			\end{tikzpicture}
		}
	\end{center}
	\caption{On the left: A junction with two ingoing and one outgoing road (2-1 node). On the right: A junction with one incoming and two outgoing roads (1-2 node).}
	\label{fig:junction}
\end{figure}
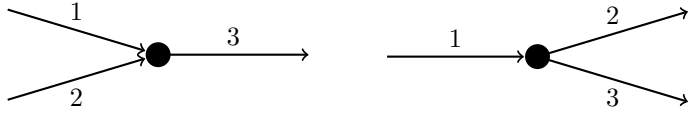

In this section we propose general coupling conditions for the relaxation  model \eqref{macro0} for different physical situations. As long as appropriate we will proceed in a similiar way as in \cite{GP06,HR} for the  definition of the coupling conditions for the ARZ-model. At the same time we state classical coupling conditions for these situations for the LWR equations.
In section \ref{macroscopiccc} the coupling conditions for the relaxation system will be related via  a matched asymptotic expansion   procedure  to the  coupling conditions 
for the LWR model on networks. There, we consider  the case of a merging junction with a special merge condition analytically. The other cases will be investigated numerically in section \ref{Numerical results}.
 We restrict ourselves here to the case of   junctions with either two ingoing and one outgoing lane (merging junction) or a junction with two outgoing and one ingoing lane (diverging junction), see  Figure \ref{fig:junction}.

As  on each road, there is exactly one outgoing characteristic family for the relaxation system, we have to provide  three conditions  at a junction connecting three roads. 
We denote the quantities at the junctions by $\rho^i,q^i, i=1,2,3$ and the corresponding Riemann invariants by $z^i,w^i$.
In any case the conservation of mass will be imposed, i.e. all cars entering a junction via one of the incoming roads
will exit on the outgoing road.
The other two conditions depend on the physical situation under consideration. We consider first the case of a 2-1 node.

\subsection{Merging lanes}

In this case the Riemann invariants $z^1,z^2,w^3$ are prescribed  and the admissible states at the junction  in $(\rho,q)$-plane fulfill

\begin{align}
\label{eq:merge_char}
q^1 = z^1 (1-\rho^1)\\
q^2=z^2 (1-\rho^2)  \nonumber\\
q^3 = \rho^3 -w^3 \nonumber.
\end{align}

The coupling conditions are the balance of fluxes
\begin{align}
\label{eq:merge_1}
q^1+q^2=q^3
\end{align}

and  a relation for the Riemann invariants for the first characteristic family  $z$ (which is also  the momentum flux of the conservative variable).  For merging junctions, one might simply  choose the balance of $z$, i.e. $z^3=z^1+z^2$ or more general

\begin{align}
\label{eq:merge_2}
z^3 = g(z^1, z^2 , w^3)
\end{align}
with some function $g$. We use as a simple example $g= z^1+z^2$ for the investigations in section \ref{macroscopiccc} and \ref{Numerical results}
and also for the illustrations in Figures \ref{fairRiemann} and \ref{state1}.
\begin{remark}
Compare \cite{HR} for the corresponding balance condition of the momentum flux for the ARZ-equations
$$
q^3z_R^3 = q^1z_R^1 + q^2z_R^2.
$$
Note that this condition is required to obtain weak solutions on the network in the sense of \cite{HR95}.
The condition $z^3=z^1+z^2$ is the analogue to this condition for the equations considered here.
\end{remark}

To continue, one more relation is
needed to solve the Riemann problem at the junction uniquely. We consider two approaches.
The first   approach assumes a   relation of the incoming fluxes  given by the prescribed quantities $z^1,z^2,w^3$.
That means, additionally to \ref{eq:merge_1}  and \ref{eq:merge_2} we consider  the general  condition
\begin{align}
\label{eq:merge_3}
\frac{q^1}{q^2} = f  
\end{align}
with 
$f = f(z^1,z^2,w^3)$
or equivalently 
\begin{align}
\label{fluxes}
q^{1} = \frac{f}{1+f} q^3\\
q^2= \frac{1}{1+f} q^3. \nonumber
\end{align}
As an example, we use simply 
\begin{align}
\label{eqz}
f= \frac{z^1}{z^2},
\end{align}
that means, that  the relation between $q^1$ and $q^2$ is given by the relation of the corresponding Riemann invariants.

The second, more general approach describes a lane merging via a condition on the flux $q^1$.
We prescribe a relation
\begin{align}
\label{fluxesgeneral}
q^{1} = F(q^3;z^1,z^2,w^3) \le \rho^1.
\end{align}
For example, for a merging with a priority lane, here lane 1,  we consider 
 additionally to the conditions \ref{eq:merge_1}  and \ref{eq:merge_2},
the following condition on the  flux  $q^1$ 
\begin{align}
\label{eq:merge_4}
q^1 = \min\{q^3,\rho^1\} .
\end{align}
This can be rewritten as 
\begin{align}
\label{priore}
q^1 = \bar q^1 
\end{align}
with $\bar q^1 = \min(q^3,q^{max}(z^1))$, where the maximal possible flux for the Riemann invariant $z$ is defined as $q^{max}(z)= \frac{z}{1+z}$.

	A more general condition giving a partial priority to one of the roads depending on a value  $P \in [0,1]$ is given by
	a convex combination of the priority conditions for lane 1 and 2, i.e.
	\begin{align}
	\label{P}
	q^1= (1-P) \bar q^1+ P  (q^3-\bar q^2)
	\end{align}
	See e.g. \cite{KCGG}
	for a partial priority merging condition for the ARZ equations.
	
\begin{figure}[h]
	\center
	\externaltikz{proportion}{
		\begin{tikzpicture}[scale=0.65]
		\def\rhoa{0.2}
		\def\rhob{0.3}
		\def\rhoc{0.6}
		
		\pgfmathsetmacro{\za}{\rhoa}
		\pgfmathsetmacro{\zb}{\rhob}
		\pgfmathsetmacro{\zc}{\za+\zb}
		\pgfmathsetmacro{\wc}{\rhoc*\rhoc}
		\pgfmathsetmacro{\barrho}{(\wc+\zc)/(1+\zc)}

		\pgfmathsetmacro{\frhoc}{\rhoc*(1-\rhoc)}
		
		\pgfmathsetmacro{\barrhozero}{(1+sqrt(1-2*\frhoc))/2}
		
		\begin{axis}[
		legend style = {at={(1,1)}, xshift=-0.1cm, yshift=0.1cm, anchor=south east},
		legend columns= 2,	
		xlabel = ${z^2}$,		
		ylabel = ${q^1 / q^3}$,
		]
		\addplot[color = blue!0!red,thick] file{Data/case1.txt};
		\addlegendentry{$f  = z^1/z^2$}
		\addplot[color = blue!33!red,thick] file{Data/case2.txt};
		\addlegendentry{$P=1/2$}
		\end{axis}
		\end{tikzpicture}
	}
	\caption{Proportion  $\frac{q^1}{q3}$ of the flux in lane 1 for $z^2 \in [0,0.5]$ and $z^1=0.4$. Coupling conditions 
	(\ref{eqz}) and (\ref{mixed2}) are shown.}
	\label{alpha}
	\centering
\end{figure}
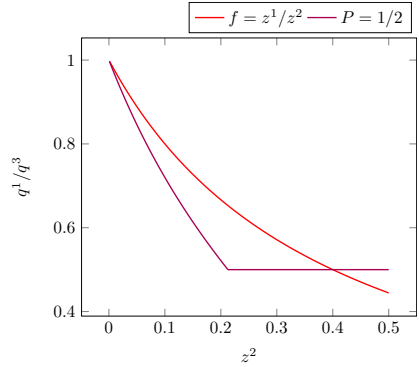

	This leads to  a distribution depending on  the maximal possible fluxes on the two ingoing roads as long as the maximal flux 
	given by $q^3$ is not exceeded and  to a distribution acording to the priority value $P$,  if both maximal possible fluxes are larger than $q^3$.
	For the symmetric case $P=\frac{1}{2}$ we obtain
	\begin{align}
	\label{mixed2}
	q^1=\frac{1}{2} \left(q^3+   \bar q^1-\bar q^2\right)	
	\end{align}
	which is easily seen to be equivalent to 	
	\begin{align}
	\label{mixed3}
	q^1 = \min \left(\bar q^1, q^3 -\min(\bar q^1, \bar q^2,\frac{q^3}{2})\right).
	\end{align}
	See Figure \ref{alpha} for the proportion $\frac{q^1}{q3}$ of the flux in  lane 1  for $z^2 $ ranging from  $0$ to $0.5$.

All these  conditions lead to a well-posed Riemann problem at the junction. This can be easily seen by solving the problems explicitly. 

\subsubsection{Solution of the Riemann problems at the junction}

\begin{figure}[h]
	\center
	\externaltikz{Statespacelayer1}{
		\begin{tikzpicture}[scale = 4]
		\def\dr{0.2}
		\node[below] at (0,0) {$0$};
		\node[left] at (0,0) {$0$};
		\node[below] at (1,0) {$1$};
		\node[left] at (0,1) {$1$};
		\draw(1,0)--(1,1);
		\draw(0,0)--(1,1);
		\draw(-0.02,1)--(0.02,1);
		
		\draw(0.68,0)--(0.68,0.48);
		\node[below] at (0.68,0) {$\bar \rho$};
		
		\draw[red](0.2,0.0)--(1,0.8);
		
		\draw(0.48,0.63)--(0.63,0.55);	
		\node[below] at (0.42,0.76) {$z^3$};	
		
		\draw[red](1.0,0.0)--(0.485,0.485);
		
		\draw(0.2,0.32)--(0.45,0.31);	
		\node[below] at (0.16,0.45) {$z^2$};		
		
		\draw[red](1.0,0.0)--(0.36,0.36);	
		
		\draw(0.3,0.47)--(0.51,0.46);	
		\node[below] at (0.27,0.6) {$z^1$};		
		
		\draw[red](1.0,0.0)--(0.6,0.6);		
		
		\draw(0.63,0.78)--(0.75,0.55);	
		\node[below] at (0.58,0.9) {$w^3$};

		\draw[](0.68,0.48)--(1.1,0.48);
		
		\node[right] at (1.1,0.48) {$q^3$};	
		
		\draw[](0.68,0.3)--(1.1,0.3);
		
		\node[right] at (1.1,0.3) {$q^1$};		
		
		\draw[](0.68,0.18)--(1.1,0.18);	
		
		\node[right] at (1.1,0.18) {$q^2$};		
		\draw[->](0,0)--(1.2,0) node[below]{$\rho$};
		\draw[->](0,0)--(0,1.2) node[left]{$q$};
		\end{tikzpicture}
		
	}
	\caption{Solution of Riemann problems for  merging with $f=z^1/z^2$.}
	\label{fairRiemann}
\end{figure}
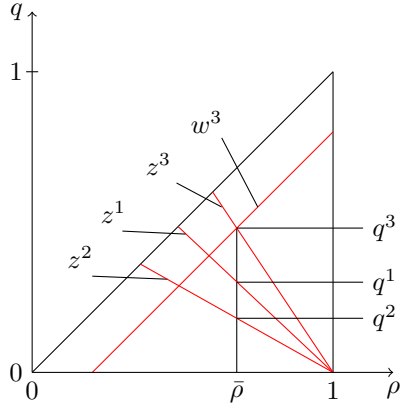

Condition \ref{eq:merge_2} directly gives   $(\rho^3,q^3)$, since 
\begin{align*}
\rho^3 -w^3 = z^3(1-\rho^3)
\end{align*}
and
\begin{align}
\rho^3  = \frac{w^3+z^3}{1+z^3} 
\end{align}
and 
\begin{align}
q^3  = z^3 (1-\rho^3 )= \frac{z^3}{1+z^3} (1-w^3).
\end{align}
See  Figure \ref{fairRiemann}. for an illustration using  the relation  (\ref{eqz}),

We obtain from \ref{fluxes}  the values of $q^1$ and $q^2$ and then  $\rho^1$ and  $ \rho^2$ from 
the characteristic equations.

Note that using the relation  (\ref{eqz}) simply gives 
that $\rho^1 = \rho^2$ and then from \ref{eq:merge_2} one obtains the equality of densities
$\rho^3 = \rho^1 = \rho^2=\bar \rho$.

For the second approach we obtain
$\rho^3,q^3$ as before. Then $q^1$ from (\ref{fluxesgeneral}) and then directly the remaining quantities.

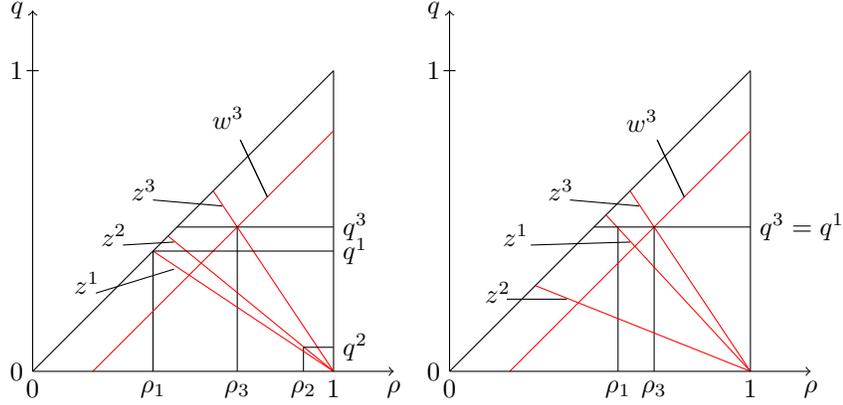
\begin{figure}[h]
	\center
	\externaltikz{Statespacelayer2}{	
		\begin{tikzpicture}[scale = 4]
		\def\dr{0.2}
		\node[below] at (0,0) {$0$};
		\node[left] at (0,0) {$0$};
		\node[below] at (1,0) {$1$};
		\node[left] at (0,1) {$1$};
		\draw(1,0)--(1,1);
		\draw(0,0)--(1,1);
		\draw(-0.02,1)--(0.02,1);
		
		\draw(0.68,0)--(0.68,0.48);
		\node[below] at (0.68,0) {$ \rho_3$};
		
		\draw(0.4,0)--(0.4,0.4);
		\node[below] at (0.4,0) {$ \rho_1$};
		
		\draw(0.9,0)--(0.9,0.08);
		\node[below] at (0.9,0) {$ \rho_2$};

		\draw[red](0.2,0.0)--(1,0.8);
		
		\draw(0.78,0.58)--(0.69,0.77);	
		\node[below] at (0.37,0.67) {$z^3$};		
		\draw[red](1.0,0.0)--(0.6,0.6);

		\draw(0.21,0.26)--(0.47,0.34);	
		\node[below] at (0.18,0.36) {$z^1$};			
		\draw[red](1.0,0.0)--(0.45,0.45);	
		
		\draw(0.3,0.42)--(0.47,0.43);	
		\node[below] at (0.27,0.52) {$z^2$};				
		\draw[red](1.0,0.0)--(0.4,0.4);		
		
		\draw(0.4,0.57)--(0.63,0.55);	
		\node[below] at (0.65,0.91) {$w^3$};

		\draw[](0.48,0.48)--(1.,0.48);
		
		\draw(0.78,0.58)--(0.69,0.77);	
		\node[right] at (1.,0.48) {$q^3$};	
		
		\draw[](0.4,0.4)--(1.0,0.4);
		\node[right] at (1.0,0.4) {$q^1$};

		\draw[](0.9,0.08)--(1.,0.08);
		\node[right] at (1.,0.08) {$q^2$};		
		\draw[->](0,0)--(1.2,0) node[below]{$\rho$};
		\draw[->](0,0)--(0,1.2) node[left]{$q$};
		\end{tikzpicture}
		\begin{tikzpicture}[scale = 4]
		\def\dr{0.2}
		\node[below] at (0,0) {$0$};
		\node[left] at (0,0) {$0$};
		\node[below] at (1,0) {$1$};
		\node[left] at (0,1) {$1$};
		\draw(1,0)--(1,1);
		\draw(0,0)--(1,1);
		\draw(-0.02,1)--(0.02,1);
		
		\draw(0.68,0)--(0.68,0.48);
		\node[below] at (0.68,0) {$ \rho_3$};
		
		\draw(0.56,0)--(0.56,0.48);
		\node[below] at (0.56,0) {$ \rho_1$};

		\draw[red](0.2,0.0)--(1,0.8);
		
		\draw(0.78,0.58)--(0.69,0.77);	
		\node[below] at (0.37,0.67) {$z^3$};		
		\draw[red](1.0,0.0)--(0.6,0.6);

		\draw(0.27,0.42)--(0.6,0.43);	
		
		\node[below] at (0.22,0.52) {$z^1$};			
		\draw[red](1.0,0.0)--(0.52,0.52);

		\draw(0.19,0.24)--(0.39,0.24);
		\node[below] at (0.16,0.33) {$z^2$};				
		\draw[red](1.0,0.0)--(0.285,0.285);		
		
		\draw(0.4,0.57)--(0.63,0.55);	
		\node[below] at (0.64,0.9) {$w^3$};

		\draw[](0.48,0.48)--(1.,0.48);
		
		\draw(0.78,0.58)--(0.69,0.77);	
		\node[right] at (1.,0.48) {$q^3=q^1$};

		\draw[->](0,0)--(1.2,0) node[below]{$\rho$};
		\draw[->](0,0)--(0,1.2) node[left]{$q$};
		\end{tikzpicture}

	}
	\caption{Solution of Riemann problems for junction with  priority lane. On the left $\frac{z^1}{1+z^1} \le q^3$. On the right $\frac{z^1}{1+z^1} \ge q^3$.}
	\label{state1}
\end{figure}

\begin{example}
For the priority  condition  (\ref{eq:merge_4}) we have to distinguish two cases. Either $\frac{z^1}{1+z^1} \ge q^3$
then $q^1=\min\{q^3,\rho^1\}=q^3$. This leads to $q^2=0$ and $\rho^2=1$. 
If $\frac{z^1}{1+z^1} \le q^3$
then $q^1=\min\{q^3,\rho^1\}=\rho^1$ and we obtain $\rho^1 = \frac{z^1}{1+z^1}$ and $q^2= q^3-q^1$.  See Figure \ref{state1} for  condition (\ref{eq:merge_4}).
\end{example}

\subsubsection{LWR-conditions}
\label{LWRcond}

\begin{figure}[h]
	\center
	\externaltikz{supplydemand}{
		\begin{tikzpicture}[scale = 3]
		\def \rhobar {0.3}	
		\def \rhostar {0.5}
		\draw[->] (0,0)--(1.2,0) node[below]{$\rho$};
		\draw[->] (0,0)--(0,1.2) node[left]{$F(\rho)$};
		\draw[dashed] (1,1.2)--(1,0.0) node[below]{$1$};
		\draw[black,line width=1pt,domain=\rhostar:1,smooth,variable=\x,] plot ({\x},{1}) ;
		\draw[black,line width=1pt,domain=0.0:\rhostar,smooth,variable=\x,] plot ({\x},{4*\x*(1-\x)}) ;
		\draw[dashed,black,line width=1pt,domain=\rhostar:1,smooth,variable=\x,] plot ({\x},{4*\x*(1-\x)}) ;
		\draw[dashed] (\rhostar,{4*\rhostar*(1-\rhostar)})--(\rhostar,0) node[below]{$\rho^*$}node at (0.5,1.1) {$\sigma$};
		\end{tikzpicture}
		\hspace{0.5cm}
		\begin{tikzpicture}[scale = 3]
		\def \rhobar {0.3}	
		\def \rhostar {0.5}
		\draw[->] (0,0)--(1.2,0) node[below]{$\rho$};
		\draw[->] (0,0)--(0,1.2) node[left]{$F(\rho)$};
		\draw[dashed] (1,1.2)--(1,0.0) node[below]{$1$};
		\draw[black,line width=1pt,domain=0.0:\rhostar,smooth,variable=\x,] plot ({\x},{1}) ;
		\draw[black,line width=1pt,domain=\rhostar:1,smooth,variable=\x,] plot ({\x},{4*\x*(1-\x)}) ;
		\draw[dashed,black,line width=1pt,domain=0.0:\rhostar,smooth,variable=\x,] plot ({\x},{4*\x*(1-\x)}) ;
		\draw[dashed] (\rhostar,{4*\rhostar*(1-\rhostar)})--(\rhostar,0) node[below]{$\rho^*$}node at (0.5,1.1) {$\sigma$};;
		\end{tikzpicture}
		
	}
	
	\caption{Supply- and demand functions $c^i$ for ingoing (left) and outgoing (right) roads.}
	\label{fig:supply}
\end{figure}
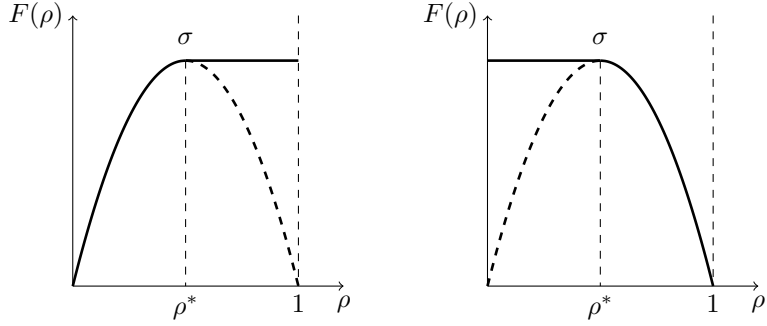

The classical conditions for LWR networks describing 
'fair merging' with equal priority and  situations with a priority lane are given by the following.
We use the supply-demand representation \cite{L} and denote the sets of valid resulting LWR fluxes $C^i$ by $\Omega^i$, compare \cite{CGP05,L,Dag1,Dag2,HR} and Figure \ref{fig:supply}. 
For the incoming roads $i=1,2$ this is
\begin{align*}
\rho_B^i \le \rho^\star \Rightarrow \Omega^i = [0, F(\rho_B^i)] &&\text{ and }&&
\rho_B^i \ge \rho^\star \Rightarrow \Omega^i = [0, \sigma]\ .	 
\end{align*}
For the outgoing road $i=3$
\begin{align*}
\rho_B^i \le\rho^\star \Rightarrow \Omega^i = [0, \sigma]  &&\text{ and }&&
\rho_B^i\ge  \rho^\star \Rightarrow \Omega^i = [0, F(\rho_B^i)]\ .	 
\end{align*}
We define the maximal admissible flux $c^i$ such that $\Omega^i = [0, c^i]$.

The coupling conditions are in both cases  given by the balance of fluxes  $C^3 = C^1+C^2$. Moreover, for 
$c^1 +c^2 \le c^3$ the conditions
\begin{align*}
C^1= c^1, C^2 = c^2\ 
\end{align*}
are used. For $c^1+c^2\geq c^3 $ we have the 'fair merging' conditions 
\begin{align}
\label{fair}
\begin{aligned}
C^i =
\min\left(c^i,c^3-\min\left(c^1,c^2,\frac{c^3}{2}\right)\right),
\qquad i = 1,2\ 
\end{aligned}
\end{align}
and in the case, where  lane 1 is a  
priority lane,
\begin{align}\label{prio}
C^1 &= \min\Big(  c^1, c^3 \Big)\\
C^2 &= c^3-C^1= \max \Big(c^3-c^1,0 \Big)  \nonumber .
\end{align}

\begin{remark}
That means, in the 'fair merging' case   the merging is  symmetric, if both incoming roads have a  flux which is larger than their share  in the outgoing road. %Otherwise, only the flux in the road with larger flux is reduced.
We refer, for example, to \cite{GPBook} for such   coupling conditions for LWR networks.
\end{remark}

In sections \ref{macroscopiccc} and \ref{Numerical results} we consider the limit of the above coupling conditions for the relaxation system  as
$\epsilon \rightarrow 0$ analytically and numerically.
In section \ref{macroscopiccc} an analytical investigation of the asymptotic procedure is given for  conditions  (\ref{eqz}) and it is shown that  the LWR condition (\ref{fair})  is obtained in the limit
as $\epsilon $ goes to 0. 
 The numerical experiments in Section \ref{Numerical results} show that the macroscopic fair merging conditions (\ref{fair}) are obtained not only from conditions (\ref{eqz}), but 
 also, for example,  from  condition (\ref{mixed2}), which is less surprising comparing  (\ref{mixed3}) and  (\ref{mixed2}).  There is obviously a larger range of conditions on the  level of the relaxation system leading to the macroscopic 
 conditions (\ref{fair}).
 We note that the   coupling condition  for the relaxation system modelling a priority lane  (\ref{eq:merge_4})  leads to the  corresponding macroscopic  condition (\ref{prio}).

\subsection{Diverging  lanes}
\label{div}
In this case the quantities 
 $z^1,w^2,w^3$ are given at the junction
and the admissible states fulfill
\begin{align}
\label{eq:div_char}
q^1 = z^1 (1-\rho^1)\\
q^2=\rho^2-w^2 \nonumber\\
q^3 = \rho^3 -w^3. \nonumber
\end{align}

Again, the coupling conditions contain  in all cases the balance  of fluxes
\begin{align}
\label{eq:div_1}
q^1=q^2+q^3.
\end{align}

To ensure that the resulting values at the junction remain in the physical domain $0 \le \rho \le 1, 0 \le q \le \rho$, we choose $q^1$ as 
$$
q^1 = \min \{  \bar q, \rho^1 \}.
$$
where  $\bar q$ is determined such that a prescribed  relation 
\begin{align}
\label{eq:div_0}
z^1 =  g(z^2,z^3)
\end{align}
is fulfilled. 
\begin{example}
	As an example we use, in Figures \ref{div1} and \ref{div2} and in the numerical experiments,  the function  $g = z^2+z^3$.
	Note that in this case the balance of the momentum flux is  only fulfilled as long as this condition leads to values inside the domain
	$ 0\le \rho \le 1, 0 \le q \le \rho $.
	For  a more detailed discussion, see Section \ref{relarz}.
\end{example}

Additionally, either the relation of   the outgoing fluxes is prescribed, i.e.  
\begin{align}
\label{eq:div_2}
\frac{q^2}{q^3} = f 
\end{align}
with $f = f(w^2,w^3,z^1) \in \R_+$, which is equivalent to
\begin{align}
\label{eq:div_2b}
q^2  =  \frac{f}{1+f} q^1, \; q^3= \frac{1}{1+f} q^1.
\end{align}
Or we use  an additive  relation 
\begin{align}
\label{eq:div_3}
q^2-q^3 = f
\end{align}
with $f = f(w^2,w^3,z^1) \in [-1,1]$.
(\ref{eq:div_3})  can be rewritten as  
\begin{align}
\label{eq:div_3b}
 q^{2} = \frac{q^1}{2} +\frac{f}{2} \; \mbox{and} \; q^{3} =\frac{q^1}{2} -\frac{f}{2} .
\end{align}

\begin{remark}
		\label{typ}
		In the first case,  typically, one prescribes a fixed relation 
		\begin{align}
		\label{pref}
		f= \frac{\alpha}{1-\alpha},
		\end{align}
		where $	\alpha \in [0,1]$ is given by the drivers preferences
		to go to one of the roads.
		This  can be rewritten as 
		\begin{align}
		\label{eq:div_2ba}
		q^{2} = \alpha q^1 \; \mbox{and} \; q^{3} = (1-\alpha) q^1 .
		\end{align}

As an example for the function $f$ in the second case, we use 
\begin{align}
\label{ex2}
f= w^3-w^2.
\end{align}
Such a  relation describes  a distribution of the outgoing fluxes adapted  to the situation in the outgoing roads
without drivers preferences. 
\end{remark}

For a comparison of these conditions with the ARZ equations, as given for example   in \cite{HR}, see section
\ref{relarz}.

Again we discuss the explicit solutions of the Riemann problems at the junction.

\subsubsection{Solution of the Riemann problems}
In the first case, we determine $q^1$ as 
$$
q^1 = \min \{\frac{z^1}{1+z^1}, \bar q  \}
$$
where $\bar q$ is determined by solving (in general numerically) the equation
\begin{align}
\label{eq1}
z^1 = g\left(\frac{f \bar q}{(1-w^2)(1+f) - f \bar q },\frac{ \bar q}{(1-w^3)(1+f) - \bar q }\right),
\end{align}
which is assumed to have  a unique solution $ q \ge 0 $ in a range, such that $\rho^2 = w^2+\frac{f}{1+f} \bar q$ and $\rho_3=w^3+ \frac{1}{1+f} \bar q $ are in $[0,1)$.

Then $q^2$ and $q^3$ are obtained straightforwardly. This yields finally  $\rho^1, \rho^2$ and $\rho^3$
due to the characteristic equations.
\begin{example}
For example, for  $g = z^2+z^3$ and $f = 1$, i.e. $\alpha =\frac{1}{2}$ and $w^2= \bar w =w^3$ we have explicitly for 
$\bar w \ge \frac{z^1}{2(1+z^1)}  $
$$
q^1= \bar q = \frac{2 z^1}{2+z^1} (1- \bar w)  
$$
and $\rho^1=\rho^2=\rho^3$.
For
$\bar w \le \frac{z^1}{2(1+z^1)}  $
we have $
q^1= \rho^1= \frac{ z^1}{1+z^1}  
$ and  $\rho^2=\rho^3$. See Figure \ref{div1} for an illustration in phase-space.
Note that for  $g$ as before, but general $\alpha\in [0,1]$ and $w^2,w^3 \in [0,1]$,  equation (\ref{eq1}) is equivalent to 
the quadratic equation $$
z^1 (a- \bar q )(b-\bar q )=  \bar q \left(  (a +b) -  2\bar q \right)   
$$
with 
$a= \frac{1-w^{2}}{\alpha}, b= \frac{1-w^{3}}{1-\alpha}$. This equation is easily seen to have a unique solution in the range 
$0 \le \bar q  \le \min (a,b)$, which is equivalent to $\rho^2 = w^2+\alpha \bar q$ and $\rho^2=w^3+(1-\alpha) \bar q$ in $[0,1]$.
\end{example}

For the second relation in Section \ref{div}, i.e. for  (\ref{eq:div_3}) we have 
$$
q^1 = \min \{\frac{z^1}{1+z^1},  \bar q \},
$$
where $\bar q$ is determined from the equation
$$
z^1 =  g\left(\frac{ \bar q+f}{2(1-\rho^2)}, \frac{\bar q-f}
{2(1-\rho^3)}\right)
$$
with $\rho^2 = q^2 + w^2= \frac{\bar q + f }{2}+w^2$ and $\rho^3 = q^3 + w^3= \frac{\bar q -f }{2}+w^3$.
 Again the equation is assumed to have  a unique solution $ q  \ge 0$  such that 
 $\rho^1$ and $\rho^2$ above are in $[0,1)$.
The remaining quantities are obtained in a straightforward way.

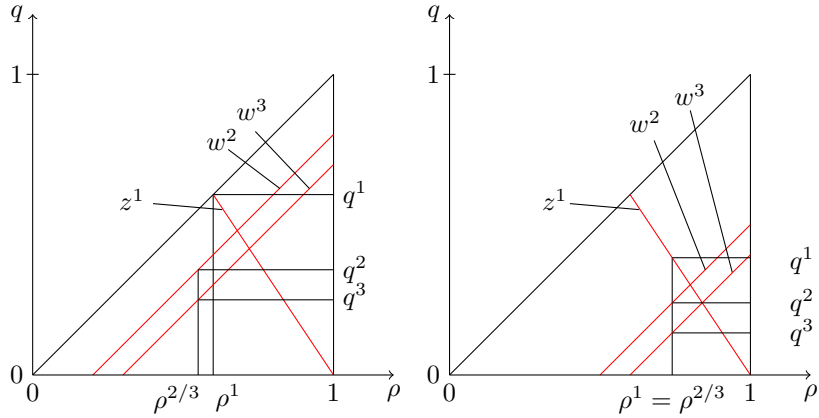
\begin{figure}[h]
	\center
	\externaltikz{Statespacelayer10}{	
		\begin{tikzpicture}[scale = 4]
		\def\dr{0.2}
		\node[below] at (0,0) {$0$};
		\node[left] at (0,0) {$0$};
		\node[below] at (1,0) {$1$};
		\node[left] at (0,1) {$1$};
		\draw(1,0)--(1,1);
		\draw(0,0)--(1,1);
		\draw(-0.02,1)--(0.02,1);
		
		\draw(0.55,0)--(0.55,0.35);
		\node[below] at (0.48,0) {$ \rho^{2/3}$};
		
		\draw(0.6,0)--(0.6,0.6);
		\node[below] at (0.65,0) {$ \rho^1$};
		
		\draw(0.4,0.57)--(0.63,0.55);		
		\node[below] at (0.33,0.65) {$z^1$};			
		\draw[red](1.0,0.0)--(0.6,0.6);	
		
		\draw[red](0.2,0.0)--(1,0.8);	
		\draw(0.65,0.75)--(0.822,0.62);
		\node[below] at (0.63,0.85) {$w^2$};		
		
		\draw[red](0.3,0.0)--(1,0.7);	
		\draw(0.74,0.81)--(0.92,0.62);	
		\node[below] at (0.73,0.95) {$w^3$};		
		
		\draw(0.55,0.35)--(1.0,0.35);	
		\node[right] at (1.,0.35) {$q^2$};

		\draw(0.55,0.25)--(1.0,0.25);	
		\node[right] at (1.,0.26) {$q^3$};

		\draw[](0.6,0.6)--(1.0,0.6);
		\node[right] at (1.0,0.6) {$q^1$};		
		
		\draw[->](0,0)--(1.2,0) node[below]{$\rho$};
		\draw[->](0,0)--(0,1.2) node[left]{$q$};
		\end{tikzpicture}
		\begin{tikzpicture}[scale = 4]
		\def\dr{0.2}
		\node[below] at (0,0) {$0$};
		\node[left] at (0,0) {$0$};
		\node[below] at (1,0) {$1$};
		\node[left] at (0,1) {$1$};
		\draw(1,0)--(1,1);
		\draw(0,0)--(1,1);
		\draw(-0.02,1)--(0.02,1);
		
		\draw(0.74,0)--(0.74,0.39);
		\node[below] at (0.74,0) {$ \rho^1=\rho^{2/3}$};
		
		\draw(0.4,0.57)--(0.63,0.55);	
		\node[below] at (0.35,0.65) {$z^1$};	
				
		\draw[red](1.0,0.0)--(0.6,0.6);	
		
		\draw[red](0.5,0.0)--(1,0.5);
		\node[below] at (0.65,0.91) {$w^{2}$};		
		
		\draw[](0.8,0.88)--(0.94,0.34);
		\draw[red](0.6,0.0)--(1,0.4);
		\node[below] at (0.8,1) {$w^{3}$};		
		
		\draw[](0.74,0.39)--(1.,0.39);
		\draw(0.68,0.78)--(0.85,0.35);	
		\node[right] at (1.1,0.38) {$q^1$};

		\draw[](0.74,0.24)--(1.,0.24);	
		\node[right] at (1.1,0.24) {$q^2$};	
		
		\draw[](0.74,0.14)--(1.,0.14);	
		\node[right] at (1.1,0.14) {$q^3$};	
		
		\draw[->](0,0)--(1.2,0) node[below]{$\rho$};
		\draw[->](0,0)--(0,1.2) node[left]{$q$};
		\end{tikzpicture}
	}
	\caption{Solution of Riemann problems for diverging junction without  drivers preferences. On the left $\frac{z^1}{1+z^1} \ge w^2+w^3$. On the right $\frac{z^1}{1+z^1} \le w^2+w^3$.}
	\label{div2}
\end{figure}

\begin{example}
Using $g = z^1+z^2$ and  $f=w^3-w^2 $ we have always $\rho^2=\rho^3 = \bar \rho$.
Moreover, 
 $\bar q$ is now determined from 
$$
z^1 = \frac{ \bar q+w^3-w^2}{2(1-\rho^2)}+\frac{\bar q-w^3+w^2}{2(1-\rho^3)}
= \frac{\bar q}{1-\frac{\bar q}{2}-\frac{w^2+w^3}{2}} .
$$
This is explicitly solved and we obtain for $w^2+w^3 \ge \frac{z^1}{1+z^1} $
$$
q^1 =  \bar q =(2-(w^2+w^3))\frac{ z^1}{2+z^1}
$$
and
$$
\rho^1 = \rho^2=\rho^3 = \frac{w^2+w^3+z^1}{2+z^1}.
$$
For $w^2+w^3 \le \frac{z^1}{1+z^1} $ we obtain 
$
q^1 = \frac{z^1}{1+z^1}$ and 
$\rho^2 = \rho^3= \frac{z^1+(1+z^1)(w^2+w^3)}{2(1+z^1)}$.
See Figure \ref{div2} for the illustration of these conditions.
\end{example}
 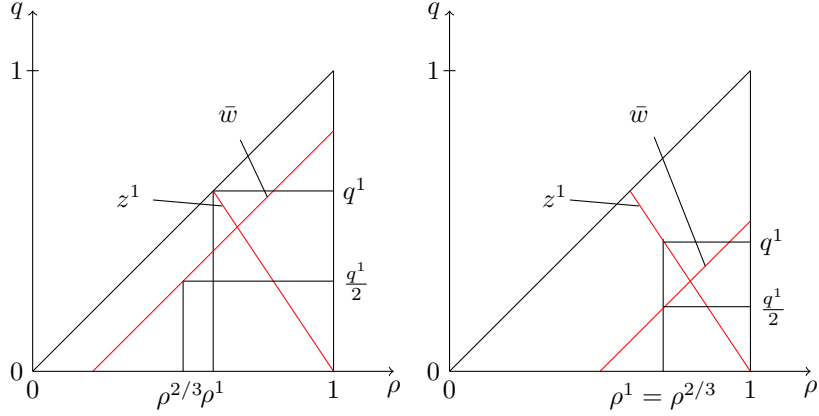
\begin{figure}[h]
	\center
	\externaltikz{Statespacelayer6}{	
		\begin{tikzpicture}[scale = 4]
		\def\dr{0.2}
		\node[below] at (0,0) {$0$};
		\node[left] at (0,0) {$0$};
		\node[below] at (1,0) {$1$};
		\node[left] at (0,1) {$1$};
		\draw(1,0)--(1,1);
		\draw(0,0)--(1,1);
		\draw(-0.02,1)--(0.02,1);

		\draw(0.6,0)--(0.6,0.6);
		\node[below] at (0.6,0) {$ \rho^1$};
		
		\draw(0.5,0)--(0.5,0.3);
		\node[below] at (0.49,0) {$ \rho^{2/3}$};

		\draw[red](0.2,0.0)--(1,0.8);

		\node[below] at (0.32,0.65) {$z^1$};			
		\draw[red](1.0,0.0)--(0.6,0.6);

		\draw(0.4,0.57)--(0.63,0.55);	
		\node[below] at (0.65,0.91) {$\bar w$};

		\draw[](0.5,0.3)--(1.,0.3);
		
		\draw(0.78,0.58)--(0.69,0.77);	
		\node[right] at (1.,0.3) {$\frac{q^1}{2}$};	
		
		\draw[](0.6,0.6)--(1.0,0.6);
		\node[right] at (1.0,0.6) {$q^1$};

		\draw[->](0,0)--(1.2,0) node[below]{$\rho$};
		\draw[->](0,0)--(0,1.2) node[left]{$q$};
		\end{tikzpicture}
		\begin{tikzpicture}[scale = 4]
		\def\dr{0.2}
		\node[below] at (0,0) {$0$};
		\node[left] at (0,0) {$0$};
		\node[below] at (1,0) {$1$};
		\node[left] at (0,1) {$1$};
		\draw(1,0)--(1,1);
		\draw(0,0)--(1,1);
		\draw(-0.02,1)--(0.02,1);
		
		\draw(0.71,0)--(0.71,0.44);
		\node[below] at (0.71,0) {$ \rho^1=\rho^{2/3}$};

		\draw[red](0.5,0.0)--(1,0.5);

		\draw(0.4,0.57)--(0.63,0.55);	
		\node[below] at (0.35,0.65) {$z^1$};			
		\draw[red](1.0,0.0)--(0.6,0.6);	
		
		\node[below] at (0.63,0.91) {$\bar w$};

		\draw[](0.71,0.43)--(1.,0.43);
		\draw(0.68,0.78)--(0.85,0.35);	
		\node[right] at (1.,0.43) {$q^1$};

		\draw[](0.71,0.215)--(1.,0.215);	
		\node[right] at (1.,0.215) {$\frac{q^1}{2}$};

		\draw[->](0,0)--(1.2,0) node[below]{$\rho$};
		\draw[->](0,0)--(0,1.2) node[left]{$q$};
		\end{tikzpicture}
	}
	\caption{Solution of Riemann problems for diverging junction with drivers preferences, $\bar w= w^2=w^3, \alpha=\frac{1}{2}$. On the left $\frac{z^1}{2(1+z^1)} \ge \bar w$. On the right $\frac{z^1}{2(1+z^1)} \le \bar w$.}
	\label{div1}
\end{figure}

\subsubsection{LWR-conditions}

Classical  coupling conditions for the LWR network for the two diverging situations are well known, see \cite{GPBook}.
Using the notation from Section \ref{LWRcond}  we always have  $C^1 = C^2+C^3$.
Additionally, we have for the situation with drivers preferences 
\begin{align}
\label{divmacro1}
C^1 = \min\Big(  c^1, \frac{1}{\alpha}  c^2 , \frac{1}{1-\alpha} c^3 \Big)
\end{align}
and 
\begin{align}
\frac{C^2}{C^3} = \frac{\alpha}{1-\alpha}.
\end{align}

Without drivers preferences
the additional conditions are for $c^2+c^3 \le c^1$
\begin{align*}
  C^2 =c^2\;,\; C^3=c^3\ 
\end{align*}
and for 
$
c^2+c^3 \ge c^1$

\begin{align}
\label{divmacro2}
C^1= c^1\;,\; C^2= \min\left(c^2,c^1-\min\left(c^2,c^3,\frac{c^1}{2}\right)\right)\ .
\end{align}
In the latter case the flow is equally distributed, if both capacities exceed the incoming flow. Otherwise, the smaller capacity 
is fully used and the lane with larger capacity has  maximal flow under these constraints.

The numerical investigations in Section \ref{Numerical results} show that condition \ref{divmacro1}  is obtained in the limit 
$\epsilon $ to $0$ from  condition (\ref{eq:div_2}) with (\ref{pref}), whereas condition \ref{divmacro2}  is obtained as the limit of condition \ref{eq:div_3} with \ref{ex2}.

\subsection{Relation to coupling conditions for the ARZ equations on networks}
\label{relarz}	
	The coupling conditions for the ARZ-equations, see 
	\cite{HR,KCGG} and many others, 
	rely on the balance of the momentum flux  $q z_R $ in the conservative formulation and the related definition of weak network solutions \cite{HR95}. 
	See the work in \cite{GP06} for an exception not requiring this condition.
	The counterpart for the present model is the balance of the momentum flux $z =\frac{q}{1-\rho}$ which has been used here for merging junctions.
	However, for diverging junctions, the balance of the quantity $z$ can in general not be prescribed anymore, if one requires that the coupling conditions lead to solutions at the nodes which remain inside the traffic domain
	$0 \le \rho \le 1$ and $0 \le q \le \rho$.
	  In general, the Riemann problem is not solvable inside this domain.
	
	We note that this is also true for the ARZ equations. In particular, the coupling conditions  for the ARZ-equations in  \cite{HR,KCGG} do not guarantee that the network solution at the nodes 
	remain in $0 \le \rho \le 1$ and $0 \le q \le \rho$. This is easily seen by  looking, for example, at the construction 
	in  \cite{HR}, where values outside this region are obtained in general.

As a final remark, we note that the present model leads to much simpler explicitly solvable conditions compared to  the ARZ model. This allows to  investigate the coupling problem in more detail, see the asymptotic investigation as $\epsilon \rightarrow 0$  in  the following sections.
	
\section{An asymptotic procedure for the relaxation system on networks in the zero relaxation limit}
\label{asyproc}
In this section we relate  the coupling conditions for the scalar conservation law  to  coupling conditions   of the nonlinear  relaxation system in the limit $\epsilon \rightarrow 0$.
This is done via a matched asymptotic expansion using  a  boundary layer analysis around the node.  
This leads to the consideration of  a  half-space problem for each lane at the node. We refer to  \cite{BSS84,BLP79,CGS,N99} for boundary layers of kinetic equations  and to \cite{AM04,LX96,NT01,WX99,WY99} for investigations of boundary layers for hyperbolic relaxation systems and kinetic equations.

The general procedure is as follows: a half space layer problem is determined by a rescaling  of the spatial coordinate on each lane  in a spatial layer 
near the node. 
The coupling conditions for the layer  problems are  given by the   coupling  condition for the relaxation system. 
Finally, the asymptotic values of the layer problems are  matched to    half Riemann problems for the macroscopic equations. 
Finally, this gives  the macroscopic coupling conditions for the LWR equations.

\subsection{The matched asymptotic expansion on the network}

We consider a single node and ingoing and outgoing lanes $[x_L^i , x_R^i]$ numbered by $i$.
The network relaxation system is given  on each lane $i$ by the relaxation equations (\ref{macro0}) for the quantities
\begin{align*}
\rho^i  (x), q^i(x)
\end{align*}
 and the coupling conditions from section \ref{sec:kineticcouplingconditions} for these values at the nodes, i.e. at
$x=x_L$ or $x=x_R$ depending, whether the lane is outgoing or ingoing.

Now, the solution of the relaxation system is approximated on each lane  by an asymptotic expansion.
For outgoing lanes this is  
\begin{align*}
\rho^i (x) \sim \rho^i_L (\frac{x-x_L}{\epsilon}) - \rho_L^i(\infty) + \rho^i_{LWR}(x)+ \mathcal{O}(\epsilon).
\end{align*}
Here  $\rho^i_L(y), y \in [0,\infty)$ is the left layer solution on lane $i$ and  $\rho^i_{LWR}$ is the LWR solution on this lane.
The  LWR value at the node  is given by 
$$
\rho_{LWR}^i = \rho_{LWR}^i (x^i_L)= \rho^i_L (\infty)= \rho^i_K.
$$
For ingoing lanes 
\begin{align*}
\rho^i (x) \sim  \rho^i_R (\frac{x_R-x}{\epsilon}) - \rho_R^i (\infty) + \rho^i_{LWR}(x)+ \mathcal{O}(\epsilon)
\end{align*}
with the layer solution $\rho^i_R(y), y \in [0,\infty)$ and 
$$
\rho_{LWR}^i = \rho_{LWR}^i (x_R^i)= \rho^i_R (\infty)= \rho^i_K.
$$
The coupling of the asymptotic expansions at the nodes means that 
$
\rho_L^i  (0), q^i_L(0)
$
for outgoing lanes and
$
\rho_R^i  (0), q^i_R(0)
$
for ingoing lanes 
fulfill the coupling conditions for the relaxation system, but not the characteristic equations (\ref{eq:merge_char}) or (\ref{eq:div_char}).

\begin{remark}
We note that the quantities $\rho_L^i  (0), q^i_L(0)$
for outgoing lanes and
$\rho_R^i  (0), q^i_R(0)$ for the ingoing lanes, which are denoted later on by $\rho_0^i, q^i_0$, are in general not equal to the values $\rho^i,q^i$ of the solution of the relaxation system at the nodes! Both fulfill the coupling conditions, but have different characteristic equations. The transition from $\rho^i,q^i$  to $\rho_R^i  (0), q^i_R(0)$ is given through a
layer in time depending  on $\epsilon$.
See Figure \ref{couplingtotal} for a graphical discussion of the situation at the node  in state space and Figure \ref{junctionvalues} for the corresponding time development of the solutions at the node.
\end{remark}

Initial conditions $(\rho^i_{init} (x) , q^i_{init} (x))$ on the lanes are  chosen in equilibrium, i.e.
$q^i_{init} (x) = F(\rho^i_{init} (x)))$ and the corresponding  values at the nodes are denoted by $\rho_B^i= \rho_{init}^i(x_L)$ or $\rho_B^i= \rho_{init}^i(x_R)$ for outgoing and ingoing roads respectively. 
In the following we investigate first the  layer equations and their asymptotic states and then the admissible half-Riemann problems  for the macroscopic equations. 

Showing the validity of the asymptotic procedure is then equivalent to  matching the asymptotic states of the 
layer problems with the admissible boundary conditions for the LWR problem  (i.e. half-Riemann problems for the LWR equations) and proving that there is a unique 
matching.

This will be done in section \ref{macroscopiccc} considering the case of a merging junction with fair-merging conditions.
The other cases will be treated numerically in Section \ref{Numerical results}. We proceed by discussing the layer problems.

\subsection{Layer solutions for the  relaxation equation}
\label{kinlayer}

We investigate  the  layers of the relaxation system at the left (outgoing lanes)  and the right (ingoing lanes) boundary.
\subsubsection{Left layer}
Consider  the left boundary of the domain being  located at $x=x_L$.
Starting from equation \eqref{macro0} and rescaling space as $y= \frac{x-x_L}{\epsilon}$ and neglecting higher order terms in $\epsilon$ one obtains 
the  layer equations  for the left boundary  for the layer solutions $(\rho_L,q_L)$ and $y \in [0, \infty)$ as 
\begin{align}
\label{layerproblem}
\begin{aligned}
\partial_y q_L &=0\\
\frac{q_L}{1-\rho_L}  \partial_y \rho_L  + (1-\frac{q_L}{1-\rho_L})\partial_y  q_L & =- \left(q_l-F(\rho_L) \right) \ .
\end{aligned}
\end{align}
This yields 
\begin{align}\label{layer}
\begin{aligned}
q_L &=C\\
\partial_y \rho_L   &= (1-\rho_L) \frac{F(\rho_L)-C}{C}\ .
\end{aligned}
\end{align}
For $0<C < F(\rho^{\star})=\sigma$, where $\rho^{\star}$ denotes the point where the maximum of $F$ is attained,
the above problem has two relevant  fix-points
$$\rho_{-} (C)\le \rho^{\star}\quad  \text{and}\quad  
\rho_+ (C)= \tau (\rho_-) \ge \rho^{\star}\ .$$
Here, $\tau(\rho)\neq \rho$ is defined by $F(\tau(\rho))= F(\rho)$.
$\rho_-$ is instable, $\rho_+ $ is stable. The domain of attraction of the stable fix point $\rho_+$  is   the interval $(\rho_-,1)$.

The third fix point   $\rho=1$ is not relevant for the further matching procedure, since  it requires $C=0$ in the macroscopic limit. In case  $C=0$ we have the instable fix point $\rho_+ = 1$ and the stable fix point $\rho_- =0$ with domain of attraction $[0,1)$. 
Moreover, we note that for $C=F(\rho^{\star})$ we have $\rho_- = \rho_+ = \rho^{\star}$ and 
all solutions with initial values above $\rho^{\star}$ converge towards $\rho^{\star}$, all other solutions  diverge.

\begin{remark}
	In case of the LWR model with $F(\rho) = \rho(1-\rho)$, see Figure \ref{figfund}, we have $$\rho_{\pm} (C) = \frac{1}{2} (1 \pm \sqrt{1-4 C})\ ,$$
	with $C< \frac{1}{4}$.
	For $C=\frac{1}{4}$ we have $\rho_- = \rho_+ = \frac{1}{2}$. Moreover, $\tau(\rho) = 1-\rho$.
	\begin{figure}[h]
		\center
		\externaltikz{LWR2}{
				\begin{tikzpicture}[scale = 3]
				\def \rhobar {0.3}	
				\def \rhostar {0.5}
				\draw[->] (0,0)--(1.2,0) node[below]{$\rho$};
				\draw[->] (0,0)--(0,1.2) node[left]{$F(\rho)$};
				\draw[dashed] (0,{4*\rhobar*(1-\rhobar)})--(1.2,{4*\rhobar*(1-\rhobar)}) node at (-0.1,{4*\rhobar*(1-\rhobar)}) {$C$}; ;
				\draw[black,line width=1pt,domain=0.0:1,smooth,variable=\x,] plot ({\x},{4*\x*(1-\x)}) ;
				\draw[dashed] (\rhobar,{4*\rhobar*(1-\rhobar)})--(\rhobar,0) node[below]{$\rho_-$};
				\draw[dashed] (1-\rhobar,{4*\rhobar*(1-\rhobar)})--(1-\rhobar,0) node[below]{$\rho_+ $};
				\draw[dashed] (\rhostar,{4*\rhostar*(1-\rhostar)})--(\rhostar,0) node[below]{$\rho^*$}node at (0.5,1.1) {$\sigma$};
				%\draw[dashed] (\rhobar,{4*\rhobar*(1-\rhobar)})--(1,0)
				%node[below]{$1$};	;
				\end{tikzpicture}
				\hspace{0.5cm}
			\begin{tikzpicture}[scale = 3]	
			\draw[->] (0,0)--(1.2,0) node[below]{$C$};
			\draw[->] (0,0)--(0,1.2) node[left]{$\rho$};
			\draw[dashed] (0,0.5)--(1.2,0.5) ;
			\draw[dashed] (1.0,0)--(1.0,1.2) ;
			 \node at (0.5,0.25) {$\rho_-$};
			  \node at (0.5,0.75) {$\rho_+$};
			 \node[left] at (0,0.5) {$\rho^\star$}; 
			 \node[below] at (1.0,0) {$\sigma$}; 
			\draw[black,line width=1pt,domain=0.0:1.0,smooth,variable=\x,samples = 100] plot ({\x},{0.5*(1-sqrt(1-\x))});
			\draw[black,line width=1pt,domain=0.0:1.0,smooth,variable=\x,samples = 100] plot ({\x},{0.5*(1+sqrt(1-\x))});
			\end{tikzpicture}
		}
		\caption{Fundamental diagram, $F(\rho)$ and fix-points of the layer problem $\rho_\mp$.}
		\label{figfund}
	\end{figure}
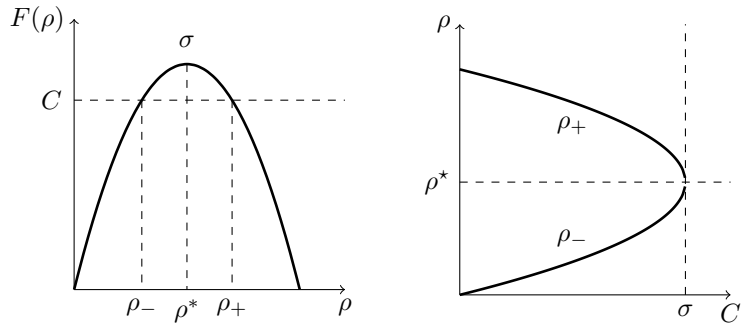
\end{remark}

\subsubsection{Right  layer}

For the right boundary at $x_R$  a scaling $y=\frac{x_R-x}{\epsilon}$ gives the layer equations for $(\rho_R,q_R)$ and  $y\in[0, \infty)$ as
\begin{align}\label{layerright}
\begin{aligned}
q_R &=C\\
- \partial_y \rho_R   &= (1-\rho_R) \frac{F(\rho_R)-C}{C}\ .
\end{aligned}
\end{align}
For  $0 < C < F(\rho^{\star})$
the above problem has again two  relevant fix points
$$\rho_{-}(C) \le \rho^{\star}\ ,\ 
\rho_+ (C) = \tau (\rho_-) \ge \rho^{\star}\; .$$
In this case 
$\rho_- $ is stable, $\rho_+ $ is instable. 
The domain of attraction of the  stable fix point $\rho_-$ is $[0,\rho_+)$.

For $C=F(\rho^{\star})=\sigma$ we have $\rho_-= \rho_+ = \rho^{\star}$ and 
all solutions with initial values below  $\rho^{\star}$ converge towards $\rho^{\star}$, all other solutions converge to not admissible states.
For $C=0$ we have the instable fix point $\rho_+ = 1$ and the stable fix point $\rho_- =0$ with domain of attraction $[0,1)$.

\subsubsection{Summary}
\label{summary}
In summary we have the following cases denoting with  $U$ the unstable fix points and with $S$ the stable ones.
Moreover, we use, as before,  the notation $\rho_K$ for the  values $\rho_L(\infty)$ and $\rho_R(\infty)$ at infinity 
and the notation $\rho_0$ for the  values at $y=0$, i.e. $\rho_L(0)$ and $\rho_R(0)$.
\\
\paragraph{Layer Problem at the left boundary} 
\begin{align*}
&\left.\begin{array}{lll}
\rho_K = \rho_-(C) \quad &\Rightarrow\quad \rho_0 =\rho_-(C), &0 \le   C< \sigma
\end{array}\right\}
&\quad \text{(U)}\\
&\left.\begin{array}{lll}
\rho_K = \rho_+(C) \quad &\Rightarrow\quad \rho_0 \in (\rho_-(C),1), &0< C < \sigma\\
\rho_K = \rho^\star \quad &\Rightarrow\quad  \rho_0 \in [\rho^\star,1),& C=\sigma\\
\rho_K = 1 \quad &\Rightarrow\quad  \rho_0 \in (0,1],& C=0
\end{array}\right\}
&\quad \text{(S)}
\end{align*}

\paragraph{The Layer Problem at the right boundary} 			
\begin{align*}
&\left.\begin{array}{lll}
\rho_K = \rho_+(C) \quad &\Rightarrow\quad \rho_0 =\rho_+(C), &0 \le  C< \sigma 
\end{array}\right\}
&\quad \text{(U)}\\
&\left.\begin{array}{lll}
\rho_K = \rho_-(C) \quad &\Rightarrow\quad \rho_0 \in [0,\rho_+(C)), &0 < C < \sigma\\
\rho_K = \rho^\star \quad &\Rightarrow\quad  \rho_0 \in [0,\rho^\star],& C=\sigma\\
\rho_K = 0 \quad &\Rightarrow\quad  \rho_0 \in [0,1),& C=0
\end{array}\right\}
&\quad \text{(S)}
\end{align*}

We use for the three cases  of the stable fix point (S) the notation
$$
\rho_K = \rho_+(C) \quad \Rightarrow\quad  \rho(0) \in \lceil \rho_-(C),1 \rfloor, 0 \le  C \le  \sigma
$$
for the left boundary and 
$$
\rho_K = \rho_-(C) \quad \Rightarrow\quad  \rho(0) \in \lceil 0, \rho_+(C) \rfloor, 0 \le  C \le  \sigma
$$
for the right boundary.

\subsection{Half-Riemann problems for the limit conservation law}
\label{Riemann}
We consider the limit conservation law  $ \partial_t \rho + \partial_x F(\rho)=0$. The initial trace at the boundary of the scalar equation is as before denoted  by $\rho_B$.
For a matching of layer solutions and solutions of the scalar conservation, those states $\rho_K$ are required which can be connected to $\rho_B$ with LWR-waves (shocks and rarefaction waves) with
non-negative velocity at the left and non-positive velocity at the right boundary. They  are summarized in the following:\\

\paragraph{The half-Riemann problem at the left boundary}

\begin{align*}
\rho_B&\leq \rho^\star \ (\text{RP 1})  \quad &\Rightarrow\quad \rho_K &\in [0,\rho^\star ]\\
\rho_B&> \rho^\star \ (\text{RP 2})  \quad &\Rightarrow\quad \rho_K &\in [0,\tau(\rho_B)]\cup\{\rho_B\}
\end{align*}		

\paragraph{The half-Riemann problem at the right boundary}		
\begin{align*}
\rho_B&\geq  \rho^\star \ (\text{RP 1})  \quad &\Rightarrow\quad \rho_K &\in [ \rho^\star ,1]\\
\rho_B&<  \rho^\star  \ (\text{RP 2})  \quad &\Rightarrow\quad \rho_K &\in \{\rho_B\}\cup [\tau(\rho_B),1]
\end{align*}

\section{The asymptotic procedure for  merging junctions}
\label{macroscopiccc}

Here we consider the merging case with coupling conditions (\ref{eqz}), i.e. the equality of densities, in detail.
First we investigate  the  layers of the relaxation system at the nodes coupled to each other via the coupling conditions and determine resulting conditions on their asymptotic states. Then, we   match these  results to Riemann solutions of the  macroscopic problems on each of the roads. The main result of this section is to show,  that the asymptotic procedure  starting from the relaxation network with the   conditions (\ref{eqz}) leads in  the limit $\epsilon \rightarrow 0$ to   the LWR-network with the  fair merging  conditions  \ref{fair}.

Assuming the boundary traces $\rho_B^i,i=1,2,3$ on the three roads to be given, we have to determine the new states $\rho_K^i$ at the node. On the one hand $\rho_K^i$ are the asymptotic states of the respective layer problems, on the other hand they are the right (for road 1 and 2) or left (for road 3) states of the half-Riemann problems for the LWR equations, i.e. the boundary conditions for the LWR equations. 
The states at the junction in the asymptotic limit (corresponding to $y=0$ for the layers) are denoted in the following  by
$\rho^i_0$ and determined together with the other values.
We have to consider eight different configurations of Riemann problems.
For  each of them all possible combinations with stable or unstable layer solutions have to be discussed.
Not admissible combinations are not listed.
 We give first 
 a detailed discussion of the coupling of the layer solutions in section  \ref{layerproof} and second a discussion of the matching of the layer solutions to the half Riemann problems  on the respective roads  in section 
\ref{proof}.

\begin{remark}
	A similiar procedure could be used for the other conditions. We limit ourselves to a numerical investigation, see Section \ref{Numerical results}.
\end{remark}

\subsection{Coupling the layers}		
\label{layerproof}

Here, we consider the coupling of the layers at the node.
The states at the junction (corresponding to the ingoing states at $y=0$ for the layers) are 
$\rho_0^i$.
Each layer can have either a stable solution (S) or an unstable solution (U).
Thus, for three lanes  we have eight possible combinations, which be denote by U/S-U/S-U/S.

\noindent{\bf Case1, U-U-U.}	
We  have $\rho_0^1 = \rho_+ (C^1), \rho_0^2= \rho_+(C^2) ,\rho_0^3= \rho_-(C^3)$, compare  subsection \ref{summary}.
The coupling conditions give
\begin{align*}
\rho_+(C^1) &= \rho_+(C^2)= \rho_-(C^3)\\
C^3&= C^1+C^2
\end{align*}	
with $0 \le  C^1,C^2,C^3< \sigma$.  

The second equality gives  $C^2= C^3=\sigma$. This is not consistent with the range of $C^2$ and $C^3$.
The case is not admissible.

\noindent{\bf Case 2, S-U-U} 
According to  subsection \ref{summary}, we have $\rho^1_0 \in \lceil 0,\rho_+ (C^1)\rfloor$ and $ \rho^2_0= \rho_+(C^2), \rho^3_0= \rho_-(C^3)$.
Inserting into the coupling conditions gives
\begin{eqnarray*}
	\rho_0^1 &= &\rho_+(C^2)= \rho_-(C^3)\\
	C^3 & =&C^1+C^2
\end{eqnarray*}	
with $0 \le  C^1 \le  \sigma$ and $0 \le  C^2,C^3< \sigma$.
Again the second  equation gives $C^2=C^3= \sigma$ which is not in the range of  $C^2,C^3$. 
The case is not admissible.

\noindent{\bf Case 3, U-S-U} 
We have $\rho_0^1 = \rho_+ (C^1), \rho_0^2 \in [0,\rho_+(C^2)), \rho_0^3= \rho_-(C^3)$.
The case is symmetric to the above and not admissible.

\noindent{\bf Case 4, U-U-S} 
We have $\rho_0^1 = \rho_+ (C^1), \rho_0^2 = \rho_+(C^2)$, $ \rho_0^3= \lceil\rho_-(C^3),1\rfloor$.
We have
\begin{align*}
\rho_+(C^1) &= \rho_+(C^2)= \rho_0^3\\
C^3&=C^1+C^2
\end{align*}
with $0 \le C^1,C^2 < \sigma$	 and $0 \le C^3 \le \sigma$.
This gives $C^1=C^2= \frac{C^3}{2}$ and $$	\rho_0^1 = 	\rho_0^2 =	\rho_0^3 =\rho_+(\frac{C^3}{2}).$$

\noindent{\bf Case 5, U-S-S} 
We have $\rho_0^1 = \rho_+ (C^1), \rho_0^2 \in \lceil0,\rho_+(C^2)\rfloor$ 
and $ \rho_0^3\in  \lceil\rho_-(C^3),1\rfloor$ with 
$0 \le C^1 < \sigma$	 and $0 \le C^2, C^3 \le \sigma$.
We have 
\begin{align*}
\rho_+(C^1) &= \rho_0^2= \rho_0^3\\
C^3&=C^1+C^2.
\end{align*}	
This gives $\rho_0^2 = \rho_0^3 = \rho_+(C^1) =\rho_+(C^3-C^2)$ with
the requirement $0 \le C^3-C^2 \le \sigma$ or
$ C^3 \ge C^2 $ and $\rho_0^2 = \rho_0^3 =\rho_+(C^3-C^2) \in [\rho_-(C^3),\rho_+(C^2)]$. It leads to $\rho_+(C^3-C^2) \le \rho_+(C^2)$ or $C^3-C^2 \ge  C_2$ or $C^3 \ge  2 C^2$.
Altogether, we have for $ 2 C^2 \le C^3$  and  $C^1 = C^3-C^2$
\begin{eqnarray*}
	\rho_0^1&=\rho_0^2 =	\rho_0^3 = \rho_+(C^3-C^2)\ .
\end{eqnarray*}	

\noindent{\bf Case 6, S-U-S} 
We have $\rho^1_0 \in  \lceil0,\rho_+ (C^1)\rfloor$  and $\rho_0^2 =\rho_+(C^2), \rho_0^3\in  \lceil\rho_-(C^3),1\rfloor$ with
$0 \le C^1 \le \sigma$	 and $0 \le C^2, C^3 < \sigma$.
The case is symmetric to case 5.
For $ 2 C^1 \le C^3  $ and $C^2 = C^3-C^1$ we have
\begin{eqnarray*}
	\rho_0^1=\rho_0^2 &=	\rho_0^3 = \rho_+(C^3-C^1)\ .
\end{eqnarray*}

\noindent{\bf Case 7, S-S-U}
We have $\rho_0^1 \in \lceil 0,\rho_+ (C^1)\rfloor $  and $\rho_0^2 \in \lceil0,\rho_+(C^2)\rfloor$ 
and $ \rho_0^3 =   \rho_-(C^3)$ with $0 \le C^1,C^2 \le \sigma$	 and $0 \le  C^3 < \sigma$.
The coupling conditions give 
\begin{align*}
\rho_0^1 &= \rho_0^2= \rho_-(C^3)\\
C^3&=C^1+C^2\ .	
\end{align*}	
This gives $\rho_0^1 = \rho_0^2 = \rho_-(C^3) $
with  the condition $0 \le C^1+C^2 < \sigma$.  Thus, for $0 \le C^1+C^2 < \sigma$ we have 
$$\rho_0^1 = \rho_0^2 = \rho_0^3 = \rho_-(C^1+C^2)\ .$$

\noindent{\bf Case 8, S-S-S} 
We have $\rho_0^1 \in \lceil 0,\rho_+ (C^1)\rfloor, \rho_0^2 \in \lceil0,\rho_+(C^2)\rfloor, \rho_0^3 \in\lceil   \rho_-(C^3),1\rfloor$
with $0 \le C^1,C^2,C^3 \le \sigma$.
The conditions are  
\begin{align*}
\rho_0^1 &= \rho_0^2= \rho_0^3\\
C^3&=C^1+C^2.	
\end{align*}	
The values of $\rho^1_0 = \rho^2_0=\rho^3$ are not uniquely determined, but  they are restricted to  the interval $  [\rho_-(C^1+C^2),\min(\rho_+(C^1),\rho_+(C^2))]$.	

These considerations yield all possible combinations of layer problems at the node. Now, they have to be matched
to the half-Riemann problems at the respective lanes.

\subsection{Matching of  Riemann problem and layer equations}
\label{proof}

Assuming  the  initial states $\rho_B^i, i=1, 2,3$ to be given, we have to determine the fluxes $C^i$ and new states $\rho_K^i$ at the node. As mentioned, on the one hand $\rho_K^i$ are the asymptotic states of the respective layer problems fulfilling the conditions in the last  section \ref{layerproof}.
On the other hand they are the left (road 1 and 2) or right  (road 3) states of the half Riemann problems which have to be connected with $\rho_B^i$.
As before, $\rho_0^i$ are the states at the junction in the limit $\epsilon \rightarrow 0$.

We  consider eight different configurations for the states $\rho_B^i$ corresponding to the possible combinations of different half Riemann problems.
For  each of them all possible combinations with stable or unstable layer solutions have to be discussed.
Not admissible combinations are not listed.
We consider the  cases  ordered in terms of 
the different possible  combinations of Riemann problems on the 3 roads using the notation 
RP1/2-1/2-1/2 for the respective combination of the half Riemann problems.\\

\noindent{\bf Case 1, RP1-1-1} $\rho_B^1 \ge \rho^\star , \rho_B^2 \ge \rho^\star , \rho_B^3 \le  \rho^\star $.
From Section  \ref{Riemann} we obtain
\begin{align*}
\rho_K^1 &\in [\rho^\star,1] :
&  (U) &\text{ or } ((S) \text{ with } C^1=\sigma) \\
\rho_K^2 &\in [\rho^\star,1]: 
& (U) &\text{ or } ((S) \text{ with } C^2=\sigma)\\
\rho_K^3 &\in [0,\rho^\star]:
& (U) &\text{ or } ((S) \text{ with } C^3=\sigma)
\end{align*}	
Then, the discussion in Section \ref{layerproof} leads to 5 different cases:
\begin{enumerate}
	%	\item[{\bf SUU}] with $C^1 = \sigma$  and a contradiction to $C^1=0$.
	%	\item[{\bf USU}] with $C^2 = \sigma$  and a contradiction to $C^2=0$.
	\item[{\bf UUS}] with $C^3 = \sigma$  and  $C^1=C^2=\frac{\sigma}{2}$ and $\rho_0^3 =\rho_+(\frac{\sigma}{2})$.
	\item[{\bf USS}] with $C^2 = C^3 = \sigma$ 
	which contradicts $C^3 \ge 2 C^2$.
	\item[{\bf SUS}] with $C^1 = C^3 = \sigma$ which contradicts  $C^3 \ge 2 C^1$.
	\item[{\bf SSU}] with $C^1 = C^2 = \sigma$ and a contradiction to $C^1+C^2 \le \sigma$. 
	\item[{\bf SSS}] with $C^1 = C^2 = C^3 = \sigma$, which gives a contradiction to the balance of fluxes.		
\end{enumerate}
This gives  $C_1=C_2=\frac{\sigma}{2}, C_3=\sigma$ and $\rho_0^i = \rho_+(\frac{\sigma}{2})$.
The values for 
$\rho_K^i $ follow directly.

%\begin{align*}
%\rho_K^1 &= \rho_+(\frac{\sigma}{2}) &\rho_K^2 &=\rho_+(\frac{\sigma}{2} ) &\rho_K^3 &=  \rho_+(\sigma)=\rho^\star\\
%\rho_0^i &= \rho_+(\frac{\sigma}{2})\ . &&&& 
%\end{align*}

\vspace{0.3cm}

\noindent{\bf Case 2, RP1-1-2} $\rho_B^1 \ge \rho^\star , \rho_B^2 \ge \rho^\star , \rho_B^3 \ge  \rho^\star $.
\begin{align*}
\rho_K^1 &\in [\rho^\star,1] :
&  (U) &\text{ or } ((S) \text{ with } C^1=\sigma) \\
\rho_K^2 &\in [\rho^\star,1]: 
& (U) &\text{ or } ((S) \text{ with } C^2=\sigma)\\		
\rho_K^3 &\in [0,\tau(\rho_B^3)] \cup \{\rho_B^3\}:
& ((U) \text{ with } C^3 \le F(\rho_B^3)) &\text{ or } ((S)
\text{ with } C^3=F(\rho_B^3))
\end{align*}	
\begin{enumerate}
	%	\item[{\bf SUU}] with $C^1 = \sigma$  and a contradiction to $C^1=0$.
	%	\item[{\bf USU}] with $C^2 = \sigma$  and a contradiction to $C^2=0$.
	\item[{\bf UUS}] with $C^3 = F(\rho_B^3)$  and  $C^1=C^2= \frac{1}{2} F(\rho_B^3)$ and $\rho_0^3 =\rho_+(\frac{1}{2}F(\rho_B^3)).$
	\item[{\bf USS}] with $C^2 =  \sigma$ 
	which contradicts $C^3 \ge 2 C^2$.
	\item[{\bf SUS}] with $C^1 =  \sigma$ which contradicts  $C^3 \ge 2 C^1$.
	\item[{\bf SSU}] with $C^1 = C^2 = \sigma$ and a contradiction to $C^1+C^2 \le \sigma$. 
	\item[{\bf SSS}] with $C^1 = C^2 = C^3 = \sigma$, which gives a contradiction to the balance of fluxes.		
\end{enumerate}
This gives $C_1=C_2=F(\frac{\rho_B^3}{2}), C_3=F(\rho_B^3)$ and $\rho_0^i =  \rho_+(\frac{1}{2}F(\rho_B^3))$.
%\begin{align*}
%\rho_K^1 &= \rho_+(\frac{1}{2}F(\rho_B^3)) &\rho_K^2 &=\rho_+(\frac{1}{2}F(\rho_B^3)) &\rho_K^3 &=  \rho_+(F(\rho_B^3))=\rho_B^3\\
%\rho_0^i &=  \rho_+(\frac{1}{2}F(\rho_B^3))\ .&&&& 
%\end{align*}

\vspace{0.3cm}

\noindent{\bf Case 3, RP1-2-1} $\rho_B^1 \ge \rho^\star , \rho_B^2 \le \rho^\star , \rho_B^3 \le  \rho^\star $. 
\begin{align*}
\rho_K^1 &\in [\rho^\star,1] :
&  (U) &\text{ or } ((S) \text{ with } C^1=\sigma) \\
\rho_K^2 &\in  [\tau(\rho_B^2),1] \cup \{\rho_B^2\}: 
& ((U) \text{ with } C^2 \le F(\rho_B^2)) &\text{ or } ((S)
\text{ with } C^2=F(\rho_B^2))	\\	
\rho_K^3 &\in [\rho^\star,1]:
& (U) &\text{ or } ((S) \text{ with } C^3=\sigma)	
\end{align*}	
\begin{enumerate}
	%	\item[{\bf SUU}] with $C^1 = \sigma$  and a contradiction to $C^1=0$.
	%	\item[{\bf USU}] with $C^2 = \sigma$  and a contradiction to $C^2=0$.
	\item[{\bf UUS}] with $C^3 = \sigma$ which gives  $C^1 = C^2= \frac{\sigma}{2}$. Moreover,  $C^2\le F(\rho_B^2)$. This is possible, if $\frac{\sigma}{2} \le F(\rho_B^2)$.
	Then, $\rho_0^3 = \rho_+(\frac{\sigma}{2})$
	\item[{\bf USS}] with $C^3 =  \sigma, C^2=F(\rho_B^2)$. 
	$C^3 \ge 2 C^2$ gives the requirement  $\frac{\sigma}{2} \ge F(\rho_B^2)$. Moreover, we have $C^1 = \sigma-F(\rho_B^2)  $ and 
	$\rho_0^2 = \rho_0^3= \rho_+(C^1)$.
	\item[{\bf SUS}] with $C^1 =  \sigma$ which contradicts  $C^3 \ge 2 C^1$.
	\item[{\bf SSU}] with $C^1 =  \sigma$ and $C^2 = F(\rho_B^2)$. This is only possible for $\rho_B^2 = 0$
	and $C^2 =0$. Then $C^3 =\sigma$ and  $\rho_0^i = \rho_-(\sigma) = \rho^\star$.
	\item[{\bf SSS}] with $C^1 =  C^3 = \sigma$ and $C^2 = F(\rho_B^2)$. This gives again $\rho_B^2=0$ and
	$C^2=0$. Then  $\rho_0^i \in [\rho_-(C^1+C^2),\min(\rho_+(C^1),\rho_+(C^2))]$ gives 
	$\rho_0^i \in [\rho_-(\sigma),\rho_+(\sigma)]$. This  leaves only $\rho_0^i =\rho^\star$.
\end{enumerate}
This gives  for $\frac{\sigma}{2} \le F(\rho_B^2) $ that   $C^1 = C^2= \frac{C^3}{2}=\frac{\sigma}{2}$ and $\rho_0^i =  \rho_+(\frac{\sigma}{2})$.
%\begin{align*}
%\rho_K^1 &= \rho_+(\frac{\sigma}{2}) &\rho_K^2 &=\rho_+(\frac{\sigma}{2}) &\rho_K^3 &=  \rho_+(\sigma)=\rho^\star\\
%\rho_0^i &=  \rho_+(\frac{\sigma}{2})\ .&&&& 
%\end{align*}

For 	$\frac{\sigma}{2} \ge F(\rho_B^2) $		one has	$C^1 = \sigma-F(\rho_B^2)  $ , $C^3 =  \sigma, C^2=F(\rho_B^2)$ and 
$\rho_0^i =  \rho_+( \sigma-F(\rho_B^2) )$.

%\begin{align*}
%\rho_K^1 &= \rho_+(\sigma-F(\rho_B^2)) &\rho_K^2 &=\rho_-(F(\rho_B^2)) =\rho_B^2&\rho_K^3 &=  \rho_+(\sigma)=\rho^\star\\
%\rho_0^i &=  \rho_+( \sigma-F(\rho_B^2) )\ .&&&& 
%\end{align*}
%

\vspace{0.3cm}

\noindent{\bf Case 4, RP2-1-1} $\rho_B^1 \le \rho^\star , \rho_B^2 \ge \rho^\star , \rho_B^3 \le  \rho^\star $. 
This case is symmetric to Case 3.

We have  for $\frac{\sigma}{2} \ge F(\rho_B^1)$ that $C^1 = F(\rho_B^1) , C^3 =\sigma$ , $C^2 = \sigma - F(\rho_B^1)$ and 
$\rho_0^i =  \rho_+(\sigma-F(\rho_B^1))$.
%\begin{align*}
%\rho_K^1 &= \rho_B^1 &\rho_K^2 &=\rho_+(\sigma-F(\rho_B^1)) &\rho_K^3 &=  \rho_+(\sigma)=\rho^\star\\
%\rho_0^i &=  \rho_+(\sigma-F(\rho_B^1))\ .&&&&
%\end{align*}

For 	$\frac{\sigma}{2} \le F(\rho_B^1) $		one has	$C^1 = C^2= \frac{C^3}{2} = \frac{\sigma}{2}$ and 
$\rho_0^i =  \rho_+(\frac{\sigma}{2})$.
%\begin{align*}
%\rho_K^1 &= \rho_+(\frac{\sigma}{2}) &\rho_K^2 &=\rho_+(\frac{\sigma}{2}) &\rho_K^3 &=  \rho_+(\sigma)=\rho^\star\\
%\rho_0^i &=  \rho_+(\frac{\sigma}{2})\ .&&&& 
%\end{align*}		

\vspace{0.3cm}

\noindent{\bf Case 5, RP1-2-2} $\rho_B^1 \ge \rho^\star , \rho_B^2 \le \rho^\star , \rho_B^3 \ge  \rho^\star .$
\begin{align*}
\rho_K^1 &\in  [\rho^\star,1] : 
& (U)  &\text{ or } ((S)
\text{ with } C^1=\sigma)	\\
\rho_K^2 &\in [\tau(\rho_B^2),1] \cup \{\rho_B^2\}:  
&  ((U)\text{ with } C^2 \le F(\rho_B^2)) &\text{ or } ((S) \text{ with } C^2= F(\rho_B^2))\\
\rho_K^3 &\in [0,1-\rho_B^3] \cup \{\rho_B^3\}: 
&((U)\text{ with } C^3 \le F(\rho_B^3)) &\text{ or } ((S) \text{ with } C^3= F(\rho_B^3))
\end{align*}	
\begin{enumerate}
	%	\item[{\bf SUU}] with $C^1 = \sigma$  and   a contradiction to $C^1=0$.
	%	\item[{\bf USU}] with $C^2 = F(\rho_B^2)$  and $C^3 \le  F(\rho_B^3)$. $C^2=0$ gives $\rho_B^2= 0$ and then 
	%	$C^3 =0$ and a contradiction.
	\item[{\bf UUS}] with $C^2 \le F(\rho_B^2) $ and $C^3= F(\rho_B^3)$. If  $ F(\rho_B^3)\le 2  F(\rho_B^2)$ then
	$C^1= C^2 = \frac{F(\rho_B^3)}{2}$ and 
	$\rho_0^3 = \rho_+(\frac{C^3}{2})$.
	\item[{\bf USS}] with $C^2 = F(\rho_B^2) , C^3 = F(\rho_B^3) $. With $C^3\ge 2 C^2$ or $F(\rho_B^3) \ge F(\rho_B^2)$ we have $C^1 =C^3-C^2$. 
	\item[{\bf SUS}] with $C^1 = \sigma , C^2 \le F(\rho_B^2), C^3 = F(\rho_B^3)$, which gives  a contradiction to  $C^3 \ge 2 C^1$.
	\item[{\bf SSU}] with $C^1  = \sigma$ and $C^2 = F(\rho_B^2),C^3 \le  F(\rho_B^3) $. This is only possible for $\rho_B^2 = 0$. Then $C^3 =\sigma, \rho_B^3 = \rho^\star$ and $\rho_0^1=\rho_-(\sigma) = \rho^\star$.
	\item[{\bf SSS}] with $C^1 = \sigma$ and $C^2 =F(\rho_B^2) , C^3 =F(\rho_B^3) $. This is only possible, if $C^2 =0$ and $\rho_B^2 =0$. This yields $C^3=\sigma$ and $\rho_B^3= \rho^\star$.
	Then  $\rho_0^i \in [\rho_-(C^1+C^2),\min(\rho_+(C^1),\rho_+(C^2))]$ gives 
	$\rho_0^i \in [\rho_-(\sigma),\rho_+(\sigma)]$,	which leaves only $\rho_0^i =\rho^\star$.
\end{enumerate}

This gives  for $F(\rho_B^3)  \le  2 F(\rho_B^2)$ that  $C_1=\frac{F(\rho_B^3)}{2}=C_2, C_3=F(\rho_B^3)$ and
$\rho_0^i =  \rho_+(\frac{F(\rho_B^3)}{2})$.
%\begin{align*}
%\rho_K^1 &= \rho_+(\frac{F(\rho_B^3)}{2}) &\rho_K^2 &=\rho_+(\frac{F(\rho_B^3)}{2} ) &\rho_K^3 &=  \rho_+(F(\rho_B^3))=\rho_B^3\\
%\rho_0^i &=  \rho_+(\frac{F(\rho_B^3)}{2})\ .&&&& 
%\end{align*}

For $F(\rho_B^3)  \ge  2 F(\rho_B^2) $	one has	$C_1=F(\rho_B^3)- F(\rho_B^2), C^2 = F(\rho_B^2), C_3=F(\rho_B^3)$ and
$\rho_0^i =  \rho_+(F(\rho_B^3)- F(\rho_B^2))$.
%\begin{align*}
%\rho_K^1 &= \rho_+(F(\rho_B^3)-F(\rho_B^2)) &\rho_K^2 &=\rho_-(F(\rho_B^2)) =\rho_B^2&\rho_K^3 &=  \rho_+(F(\rho_B^3))=\rho_B^3\\
%\rho_0^i &=  \rho_+(F(\rho_B^3)- F(\rho_B^2))\ .&&&& 
%\end{align*}

\vspace{0.3cm}

\noindent{\bf Case 6, RP2-1-2} $\rho_B^1 \le \rho^\star , \rho_B^2 \ge \rho^\star , \rho_B^3 \ge  \rho^\star $.
This case is symmetric to case 5.

We have for $F(\rho_B^3)  \le  2 F(\rho_B^1)$ that $C_1=\frac{F(\rho_B^3)}{2}=C_2, C_3=F(\rho_B^3)$ and
$\rho_0^i =  \rho_+(\frac{F(\rho_B^3)}{2})$.
%\begin{align*}
%\rho_K^1 &= \rho_+(\frac{F(\rho_B^3)}{2}) &\rho_K^2 &=\rho_+(\frac{F(\rho_B^3)}{2} ) &\rho_K^3 &=  \rho_+(F(\rho_B^3))=\rho_B^3\\
%\rho_0^i &=  \rho_+(\frac{F(\rho_B^3)}{2})\ .&&&& 
%\end{align*}

For $F(\rho_B^3 ) \ge  2 F(\rho_B^1) $	one has	$C_1=F(\rho_B^1), C^2 = F(\rho_B^3)- F(\rho_B^1), C_3=F(\rho_B^3)$
and $\rho_0^i =  \rho_+(F(\rho_B^3)- F(\rho_B^1))$.
%\begin{align*}
%\rho_K^1 &= \rho_-(F(\rho_B^1)) =\rho_B^1&\rho_K^2 &=\rho_+(F(\rho_B^3)- F(\rho_B^1)) &\rho_K^3 &=  \rho_+(F(\rho_B^3))= \rho_B^3\\
%\rho_0^i &=  \rho_+(F(\rho_B^3)- F(\rho_B^1))\ .&&&& 
%\end{align*}

\vspace{0.3cm}

\noindent{\bf Case 7, RP2-2-1}  $\rho_B^1 \le \rho^\star , \rho_B^2 \le \rho^\star , \rho_B^3 \le  \rho^\star $. 
\begin{align*}
\rho_K^1 &\in [\tau(\rho_B^1),1] \cup \{\rho_B^1\}:  
&  ((U)\text{ with } C^1 \le F(\rho_B^1)) &\text{ or } ((S) \text{ with } C^1= F(\rho_B^1))\\
\rho_K^2 &\in  [\tau(\rho_B^2),1] \cup \{\rho_B^2\}: 
&((U)\text{ with } C^2 \le F(\rho_B^2)) &\text{ or } ((S) \text{ with } C^2= F(\rho_B^2))\\
\rho_K^3 &\in [0,\rho^\star] : 
& (U)  &\text{ or } ((S)
\text{ with } C^3=\sigma)	 
\end{align*}	
\begin{enumerate}
	%	\item[{\bf SUU}] with $C^1 =   F(\rho_B^1)$, $C^2  \le  F(\rho_B^2)$. $C^1=0$ and $C^2 =C^3=\sigma$ gives $\rho_B^1 =0$ and $\rho_B^2 = \rho^\star$.
	%	\item[{\bf USU}] with $C^1 \le  F(\rho_B^1)$  and $C^2 =F(\rho_B^1)$. $C^2=0, C^1=\sigma$ gives $\rho_B^1 =\rho^\star $ and
	%	$\rho_B^2 =0$.
	\item[{\bf UUS}] with $C^1 \le F(\rho_B^1) $ and $C^2 \le  F(\rho_B^2) $.  $C^3=\sigma$ yields 
	$C^1= C^2 = \frac{\sigma}{2}$, if $F(\rho_B^1) \ge \frac{\sigma}{2}$ and $F(\rho_B^2) \ge \frac{\sigma}{2}$ .
	Then $\rho_0^3 = \rho_+(\frac{C^3}{2})$.
	\item[{\bf USS}] with $C^1 \le  F(\rho_B^1) , C^3 =\sigma, C^2 =  F(\rho_B^2) $.  $C^3\ge 2 C^2$ is equivalent to
	$\frac{\sigma}{2} \ge F(\rho_B^2)$. Moreover, $C^1 = \sigma - F(\rho_B^2) $ requires $F(\rho_B^1) +F(\rho_B^2) \ge \sigma$.
	\item[{\bf SUS}]  with $C^1 =  F(\rho_B^1), C^2 \le  F(\rho_B^2), C^3 =\sigma$.  $C^3 \ge 2 C^1$ gives $ \frac{\sigma}{2} \ge F(\rho_B^1)$,  $F(\rho_B^2) \ge \frac{\sigma}{2}$ and $C^2 =\sigma - F(\rho_B^1) \ge \frac{\sigma}{2}$.
	Moreover, $\rho_0^1= \rho_+(C^3-C^1)$.
	\item[{\bf SSU}] with $C^1  = F(\rho_B^1)$ and $C^2 =  F(\rho_B^2) $.  This  gives $F(\rho_B^1)+ F(\rho_B^2)\le \sigma$ and  $\rho_0^i =
	\rho_-(C^3) $.
	\item[{\bf SSS}] with $C^1 = F(\rho_B^1) $ and $C^2 =F(\rho_B^2) , C^3 =\sigma$. This is only possible, if $ F(\rho_B^1)+ F(\rho_B^2)= \sigma$. In this case, since $ \rho_0^i\in  [\rho_-(C^1+C^2),\min(\rho_+(C^1),\rho_+(C^2))]$ we obtain
	$ \rho_0^i\in  [\rho_-(\sigma),\min(\rho_+(F(\rho_B^1)),\rho_+(F(\rho_B^2)))]$. This gives the restriction
	$ \rho_0^i\in $$ [\rho^\star,\min(\tau(\rho_B^1),\tau(\rho_B^2))]$ according to the range of $\rho_B^1, \rho_B^2 $.
\end{enumerate}
We obtain for $F(\rho_B^1)+F(\rho_B^2) \le \sigma (SSU)$ that $C_1=F(\rho_B^1), C_2=F(\rho_B^2), C_3=F(\rho_B^1)+F(\rho_B^2)$ and $\rho_0^i =  \rho_-(F(\rho_B^1)+ F(\rho_B^2))$.
%\begin{align*}
%\rho_K^1 &= \rho_-(F(\rho_B^1) )=\rho_B^1&\rho_K^2 &=\rho_-(F(\rho_B^2)) =\rho_B^2&\rho_K^3 &=  \rho_-(F(\rho_B^1)+ F(\rho_B^2))\\
%\rho_0^i &=  \rho_-(F(\rho_B^1)+ F(\rho_B^2)).&&&& \ .
%\end{align*}

For $F(\rho_B^1)+F(\rho_B^2) \ge \sigma ,F(\rho_B^1) \ge \frac{\sigma}{2}, F(\rho_B^2) \ge \frac{\sigma}{2} (UUS)$	one has	$C_1= \frac{\sigma}{2}=C^2, C^3 =\sigma $ and $\rho_0^i =  \rho_+(\frac{\sigma}{2})$.
%\begin{align*}
%\rho_K^1 &= \rho_+(\frac{\sigma}{2}) &\rho_K^2 &=\rho_+(\frac{\sigma}{2})&\rho_K^3 &=  \rho^\star\\
%\rho_0^i &=  \rho_+(\frac{\sigma}{2}).&&&& \ .
%\end{align*}

For $F(\rho_B^1)+F(\rho_B^2) \ge \sigma ,F(\rho_B^1) \le \frac{\sigma}{2}, F(\rho_B^2) \ge \frac{\sigma}{2} (SUS)$	one has	$C_1=F(\rho_B^1), C_2=\sigma-F(\rho_B^1)$, $C_3=\sigma$ and $\rho_0^i =  \rho_+(\sigma-F(\rho_B^1))$.
%\begin{align*}
%\rho_K^1 &= \rho_-(F(\rho_B^1)) =\rho_B^1&\rho_K^2 &=\rho_+(\sigma-F(\rho_B^1))&\rho_K^3 &=  \rho_+(\sigma)=\rho^\star\\
%\rho_0^i &=  \rho_+(\sigma-F(\rho_B^1)).&&&& \ .
%\end{align*}

For $F(\rho_B^1)+F(\rho_B^2) \ge \sigma ,F(\rho_B^1) \ge \frac{\sigma}{2}, F(\rho_B^2) \le \frac{\sigma}{2}(USS)$	one has	$C_1=\sigma-F(\rho_B^2), C_2=F(\rho_B^2)$, $ C_3=\sigma$ and $\rho_0^i =  \rho_+(\sigma- F(\rho_B^2))$.
%\begin{align*}
%\rho_K^1 &= \rho_+(\sigma-F(\rho_B^2)) &\rho_K^2 &=\rho_-( F(\rho_B^2)) = \rho_B^2&\rho_K^3 &=  \rho^\star\\
%\rho_0^i &=  \rho_+(\sigma- F(\rho_B^2)).&&&& \ .
%\end{align*}		

\begin{remark}
	We note that at the interfaces between the different conditions we obtain values  $\rho_0^i  \in   [\rho^\star,\min(\rho_+(F(\rho_B^1)),\rho_+(F(\rho_B^2))]$. This is exactly   the interval for the $\rho^i$-values in case (SSS).
\end{remark}
\vspace{0.3cm}

\noindent{\bf Case 8, RP2-2-2} $\rho_B^1 \le \rho^\star , \rho_B^2 \le \rho^\star , \rho_B^3 \ge  \rho^\star $.
\begin{align*}
\rho_K^1 &\in [\tau(\rho_B^1),1] \cup \{\rho_B^1\}:  
&  ((U)\text{ with } C^1 \le F(\rho_B^1)) &\text{ or } ((S) \text{ with } C^1= F(\rho_B^1))\\
\rho_K^2 &\in [\tau(\rho_B^2),1] \cup \{\rho_B^2\}:  
&  ((U)\text{ with } C^2 \le F(\rho_B^2)) &\text{ or } ((S) \text{ with } C^2= F(\rho_B^2))\\
\rho_K^3 &\in [0,\tau(\rho_B^3)] \cup \{\rho_B^3\}: 
&((U)\text{ with } C^3 \le F(\rho_B^3)) &\text{ or } ((S) \text{ with } C^3= F(\rho_B^3))
\end{align*}	
\begin{enumerate}
	%	\item[{\bf SUU}] with $C^1 =   F(\rho_B^1)$, $C^2 \le   F(\rho_B^2)$. $C^2 \le   F(\rho_B^2)$.
	%	$C^1=0$ and $C^3 =\sigma$  gives $\rho_B^1 =0$ and $\rho_B^2 = \rho_B^3 = \rho^\star$.
	%	\item[{\bf USU}] with $C^1 \le  F(\rho_B^1)$  and $C^2 =F(\rho_B^1)$, $C^3 \le F(\rho_B^3)$. $C^2=0$ gives $\rho_B^2 =0$ and $\rho_B^1 = \rho_B^3 = \rho^\star$ .
	\item[{\bf UUS}] with $C^1 \le F(\rho_B^1) $,  $C^2 \le F(\rho_B^2) $ and $C^3 = F(\rho_B^3) $.  If  $ F(\rho_B^3)\le 2  F(\rho_B^1)$ and $ F(\rho_B^3)\le 2  F(\rho_B^2)$  then
	$C^1= C^2 = \frac{F(\rho_B^3)}{2}$ and 
	$\rho_0^3 = \rho_+(\frac{C^3}{2})$.
	\item[{\bf USS}] with $C^1 \le  F(\rho_B^1) , C^2 =F(\rho_B^2), C^3 = F(\rho_B^3) $. With $C^3\ge 2 C^2$ we have 
	$ F(\rho_B^3)\ge 2  F(\rho_B^2)$ and $ F(\rho_B^3)-  F(\rho_B^2) \le F(\rho_B^1)$ 
	or $ F(\rho_B^1)+ F(\rho_B^2) \ge F(\rho_B^3)$. 
	\item[{\bf SUS}] with $C^1 =  F(\rho_B^1), C^2\le  F(\rho_B^2),C^3 = F(\rho_B^3)$.  $C^3 \ge 2 C^1$ gives $ F(\rho_B^3)\ge 2  F(\rho_B^1)$ and $F(\rho_B^3)- F(\rho_B^1) \le F(\rho_B^2)$ or  $F(\rho_B^1)+ F(\rho_B^2) \ge F(\rho_B^3)$. Moreover $\rho_0^1= \rho_+(C^3-C^1)$.
	\item[{\bf SSU}] with $C^1  = F(\rho_B^1)$ and $C^2 =  F(\rho_B^2) ,C^3 \le  F(\rho_B^3) $. 
	This is only possible for $F(\rho_B^1)+F(\rho_B^2) \le F(\rho_B^3)$. Moreover, $\rho_0^i =\rho_-(C^3) $.
	\item[{\bf SSS}] with $C^1 = F(\rho_B^1) $ and $C^2 =F(\rho_B^2), C^3 =F(\rho_B^3) $. This is only possible, if $F(\rho_B^1)+F(\rho_B^2) = F(\rho_B^3)$.
	
	We obtain
	$ \rho_0^i\in  [\rho_-(F(\rho_B^3)),\min(\rho_+(F(\rho_B^1)),\rho_+(F(\rho_B^2)))]$. This gives 
	according to the range of $\rho_B^2, \rho_B^3 $, that $ \rho_0^i\in  [\tau(\rho_B^3),\min(\tau(\rho_B^1),\tau(\rho_B^2))]$.
\end{enumerate}
This gives  for $ F(\rho_B^3)\le 2  F(\rho_B^1)$ and $ F(\rho_B^3)\le 2  F(\rho_B^2) (UUS)$
that $C_1=\frac{F(\rho_B^3)}{2}=C^2,  C_3=F(\rho_B^3)$ and $\rho_0^i =  \rho_+(\frac{F(\rho_B^3)}{2})$.
%\begin{align*}
%\rho_K^1 &= \rho_+(\frac{F(\rho_B^3)}{2}) &\rho_K^2 &=\rho_+(\frac{F(\rho_B^3)}{2} ) &\rho_K^3 &=  \rho_+(F(\rho_B^3))=\rho_B^3\\
%\rho_0^i &=  \rho_+(\frac{F(\rho_B^3)}{2})\ .&&&& 
%\end{align*}

For $ F(\rho_B^3)\ge 2  F(\rho_B^2)$ and  $ F(\rho_B^1)+ F(\rho_B^2) \ge F(\rho_B^3)(USS)$	one has	
$C^2=F(\rho_B^2)$, $C^3=F(\rho_B^3)$, $C_1=F(\rho_B^3)-F(\rho_B^2)$ and $\rho_0^i =  \rho_+(F(\rho_B^3)- F(\rho_B^2))$.
%\begin{align*}
%\rho_K^1 &= \rho_+(F(\rho_B^3)- F(\rho_B^2)) &\rho_K^2 &=\rho_-(F(\rho_B^2))&\rho_K^3 &=  \rho_+(F(\rho_B^3))\\
%\rho_0^i &=  \rho_+(F(\rho_B^3)- F(\rho_B^2))\ .&&&& 
%\end{align*}

For $F(\rho_B^3)\ge 2  F(\rho_B^1)$ and  $F(\rho_B^1)+ F(\rho_B^2) \ge F(\rho_B^3) (SUS)$	one has	
$C_1=F(\rho_B^1)$,  $C^3=F(\rho_B^3)$, $C^2=F(\rho_B^3)-F(\rho_B^1)$ and $\rho_0^i =  \rho_+(F(\rho_B^3)- F(\rho_B^1))$.
%\begin{align*}
%\rho_K^1 &= \rho_-(F(\rho_B^1))= \rho_B^1&\rho_K^2 &=\rho_+(F(\rho_B^3)- F(\rho_B^1)) &\rho_K^3 &=  \rho_+(F(\rho_B^3))=\rho_B^3\\
%\rho_0^i &=  \rho_+(F(\rho_B^3)- F(\rho_B^1))\ .&&&& 
%\end{align*}						

For  $F(\rho_B^1)+F(\rho_B^2) \le F(\rho_B^3)(SSU)$	one has	$C_1=F(\rho_B^1),C^2=F(\rho_B^2)$, $C^3=F(\rho_B^1)+F(\rho_B^2)$
and $\rho_0^i =  \rho_-(F(\rho_B^1)+F(\rho_B^2))$.
%\begin{align*}
%\rho_K^1 &= \rho_-(F(\rho_B^1)) =\rho_B^1&\rho_K^2 &=\rho_-(F(\rho_B^2))=\rho_B^2 &\rho_K^3 &=  \rho_-(F(\rho_B^1)+F(\rho_B^2))\\
%\rho_0^i &=  \rho_-(F(\rho_B^1)+F(\rho_B^2))\ .&&&& 
%\end{align*}	

%Alltogether this gives exactly the 8 cases in statement \ref{}.

\begin{remark}
	Note that the sub-cases in Case 8  partition uniquely the range of admissible states since for $0 \le x,y,z \le 1$ either ($x+y\le z$) or ($x+y\ge z$ and $z\ge2y$) or ($x+y\ge z$ and $z\ge 2x$) or ($z\le 2x$ and $z\le 2y$).	
	Moreover, note that at the interfaces between the different conditions we obtain that $\rho_0^i \in [\rho_-(F(\rho_B^3)),
	\min(\rho_+(F(\rho_B^1)),\rho_+(F(\rho_B^2)))]$. This is exactly   the interval for the $\rho^i$-values in case (SSS). 
\end{remark}

The above computations show that there is a unique matching of layer solutions and LWR solutions and that 
the asymptotic expansion lead to well defined conditions for the LWR network.

Considering only the fluxes and neglecting the information on the $\rho_0^i$ the above  result can be rewritten in a more convenient way using the supply- and demand formulation, see section \ref{LWRcond} or \cite{L}.
We obtain 

{\bf Case 1, RP1-1-1.} 
This is a case with  $c^1,c^2\le \frac{c^3}{2}: C^1 =c^1, C^2 = c^2.$

{\bf Case 2, RP1-1-2.}  
This is a case with $c^1,c^2\ge \frac{c^3}{2}: C^1 =C^2 = \frac{c^3}{2}.$

{\bf Case 3, RP1-2-1} 
We have $c^1 \ge c^2$ and  two cases:
\begin{align*}
c^2 \ge \frac{c^3}{2}&: C^1= C^2 = \frac{c^3}{2}\\
c^2 \le \frac{c^3}{2}&:  C^1= c^3-c^2, C^2 =c^2.
\end{align*}

{\bf Case 4,	RP2-1-1} Symmetric to Case 3.
We have $c^1 \le c^2$ and  two cases:
\begin{align*}
c^1 \ge \frac{c^3}{2}&: C^1= C^2 = \frac{c^3}{2}\\
c^1 \le \frac{c^3}{2}&:  C^1= c^1, C^2 =c^3 -c^1.
\end{align*}

{\bf Case 5, RP1-2-2} 
In terms of the  $c^i $ this case is the same as Case 3.

{\bf Case 6, RP2-1-2}
This case is the same as Case 4.

{\bf Case 7,	RP2-2-1}  
We have four cases:
\begin{align*}
c^1 +c^2 \le c^3&:C^1= c^1, C^2 = c^2\\
c^1 +c^2 \ge c^3,c^1 \ge \frac{c^3}{2}, c^2\ge \frac{c^3}{2} &:C^1= C^2 =\frac{c^3}{2}\\
c^1 +c^2 \ge c^3,c^1 \ge \frac{c^3}{2}, c^2\le \frac{c^3}{2}&:C^1= c^3-c^2,C^2 =c^2\\
c^1 +c^2 \ge c^3,c^1 \le \frac{c^3}{2}, c^2\ge \frac{c^3}{2}  &:C^1= c^1, C^2 =c^3-c^1.
\end{align*}

{\bf Case 8, RP2-2-2}
We obtain the same as in Case 7.\\

%\noindent All in all, only 4 different cases are left:
%
%	{\bf Case A:}
%		\begin{align*}
%		c^1 +c^2 \le c^3:C^1= c^1, C^2 = c^2\ ,
%		\end{align*}
%	\qquad {\bf Case B:}
%		\begin{align*}
%		c^1 +c^2 \ge c^3,c^1 \ge \frac{c^3}{2}, c^2\ge \frac{c^3}{2}:  C^1= C^2 =\frac{c^3}{2}\ ,
%		\end{align*}
%	\qquad {\bf Case C:}
%		\begin{align*}
%		c^1 +c^2 \ge c^3,c^1 \ge \frac{c^3}{2}, c^2\le \frac{c^3}{2}:  C^1= c^3-c^2,C^2 =c^2\ ,
%		\end{align*}
%	\qquad {\bf Case D:}
%		\begin{align*}
%		c^1 +c^2 \ge c^3,c^1 \le \frac{c^3}{2}, c^2\ge \frac{c^3}{2}:  C^1= c^1, C^2 =c^3-c^1\ .
%		\end{align*}

One oberves directly, that this result can be rewritten in the 
 more compact form given in  \ref{fair} which shows that the relaxation network with conditions (\ref{eqz})
 converges  for $\epsilon  \rightarrow 0$ to the LWR network with the fair merging  conditions (\ref{fair}).

\begin{remark}
	The above derivation shows that the  classical merge condition (\ref{fair}) for the LWR-network
	can be obtained as the asymptotic limit of  condition (\ref{eqz}) for the relaxation system, that means the equality of densities.
	The  equality of densities is not fulfilled on the macroscopic level of the conservation law,
	only the balance of fluxes is common for both levels of coupling conditions.
	Moreover, we note that this is not the only coupling condition for the relaxation system  leading to (\ref{fair}). A similiar investigation leading to the same macroscopic coupling conditions could  be performed for condition
	(\ref{mixed2}), i.e. the priority condition with priority $P=\frac{1}{2}$.
\end{remark}

A graphical sketch of the different  quantitites  in state space is given  in Figure \ref{couplingtotal} for the  special example
$\rho_B^1= 0.2$, $ \rho_B^2= 0.3$ and $\rho_B^3= 0.6$. 
One observes the difference at the junction between the values $\bar \rho= \rho^i$  found by solving the coupling conditions for the relaxation system and  the values $\bar \rho_0= \rho^i_0 $ found by the asymptotic investigation.
We note that both values fulfill the coupling conditions, but different equations (characteristic equations for the relaxation system versus layer plus LWR-wave for the limit problem) connecting them to the $\rho_B^i$.
The time development of the value at the nodes for the relaxation system with different $\varepsilon$  is shown in Figure \ref{junctionvalues}.
One observes the evolution of the values at the junction from $\bar \rho$  to  $\bar \rho_0$.
The size of the temporal layer in Figure \ref{junctionvalues} also depends on $\epsilon$.

\begin{figure}[h]
	\center
	\externaltikz{cc1}{
		\begin{tikzpicture}[scale = 9.5]
		
		\def\rhoa{0.2}
		\def\rhob{0.3}
		\def\rhoc{0.6}
		
		\pgfmathsetmacro{\za}{\rhoa}
		\pgfmathsetmacro{\zb}{\rhob}
		\pgfmathsetmacro{\zc}{\za+\zb}
		\pgfmathsetmacro{\wc}{\rhoc*\rhoc}
		\pgfmathsetmacro{\barrho}{(\wc+\zc)/(1+\zc)}

		\pgfmathsetmacro{\zam}{\za/(1+\za)}
		\pgfmathsetmacro{\zbm}{\zb/(1+\zb)}
		\pgfmathsetmacro{\zcm}{\zc/(1+\zc)}
		
		\pgfmathsetmacro{\qa}{\za*(1-\barrho)}
		\pgfmathsetmacro{\qb}{\zb*(1-\barrho)}
		\pgfmathsetmacro{\qc}{\zc*(1-\barrho)}
		
		\pgfmathsetmacro{\frhoa}{\rhoa*(1-\rhoa)}
		\pgfmathsetmacro{\frhob}{\rhob*(1-\rhob)}
		\pgfmathsetmacro{\frhoc}{\rhoc*(1-\rhoc)}
		\pgfmathsetmacro{\frhochalf}{\frhoc/2}
		
		\pgfmathsetmacro{\rhoatau}{1-\rhoa}
		\pgfmathsetmacro{\watau}{\rhoatau*\rhoatau}
		\pgfmathsetmacro{\frhoatau}{\rhoatau*(1-\rhoatau)}
		
		\pgfmathsetmacro{\barrhozero}{(1+sqrt(1-2*\frhoc))/2}
		
		\pgfmathsetmacro{\wabar}{\barrho-\qa}
		\pgfmathsetmacro{\wbbar}{\barrho-\qb}
		
		\node[below] at (0,0) {$0$};
		\node[left] at (0,0) {$0$};
		\node[below] at (1,0) {$1$};
		\draw(1,0)--(1,0.7);
		\draw(0,0)--(0.7,0.7);
		
		\node[left] at (\zcm,\zcm) {$z^3$};		
		\node[left] at (\zbm,\zbm) {$z^2$};			
		\node[left] at (\zam,\zam) {$z^1$};		
		
		\node[right] at (1,1-\wc) {$w^3$};		
		
		\draw[dashed](\barrho,\qc)--(1.0,\qc);
		\node[right] at (1,\qc) {$q^3$};	
		
		\draw[dashed](\barrho,\qb)--(1.0,\qb);
		\node[right] at (1,\qb+0.01) {$q^2$};

		\draw[dashed](\barrho,\qa)--(1.0,\qa);
		\node[right] at (1,\qa) {$q^1$};	
		
		\draw[dashed](\barrho,0)--(\barrho,0.3) node[below] at (\barrho-0.01,0) {$\bar \rho$};
		
		\draw[dashed](\barrhozero,0)--(\barrhozero,0.3) node[below] at (\barrhozero,0) {$\bar \rho_0$};	
		
		\draw[->](0,0)--(1.2,0) node[below]{$\rho$};
		\draw[->](0,0)--(0,0.7) node[left]{$q$};

		\draw[domain=\wc:1.0,smooth,variable=\x,red] plot ({\x},{\x-\wc});
		\draw[domain=\zam:1.0,smooth,variable=\x,red] plot ({\x},{(1-\x)*\za});
		\draw[domain=\zbm:1.0,smooth,variable=\x,red] plot ({\x},{(1-\x)*\zb});	
		\draw[domain=\zcm:1.0,smooth,variable=\x,red] plot ({\x},{(1-\x)*\zc});

		\draw[domain=0.0:1,smooth,variable=\x,green] plot ({\x},{\x*(1-\x)});	
		
		\draw[dashed](\rhoa,0)--(\rhoa,0.2) node[below] at (\rhoa,0.0) {$\rho_B^1$};
		\draw[dashed](\rhob,0)--(\rhob,0.25) node[below] at (\rhob,0.0) {$ \rho_B^2$};
		\draw[dashed](\rhoc,0)--(\rhoc,0.25) node[below] at (\rhoc+0.01,0.0) {$\rho_B^3$};
		\draw[dashed](0.5,0)--(0.5,0.25) node[below] at (0.5,0.0) {$ \rho^\star$};
		
		\draw[] node at (\rhoa,\frhoa) {\textsf{x}};
		\draw[] node at (\rhob,\frhob) {\textsf{x}};
		\draw[] node at (\rhoc,\frhoc) {\textsf{x}};
		
		\draw[] node at (\barrho,\qa) {\textsf{x}};
		\draw[] node at (\barrho,\qb) {\textsf{x}};
		\draw[] node at (\barrho,\qc) {\textsf{x}};
		
		\draw[] node at (\barrhozero,\frhochalf) {\textsf{x}};
		\draw[] node at (\barrhozero,\frhochalf) {\textsf{x}};
		\draw[] node at (\barrhozero,\frhoc) {\textsf{x}};
		
		\draw[dashed](\barrhozero,\frhochalf)--(1.08, \frhochalf) node[right] at (1.07,\frhochalf) {$ \frac{F(\rho_B^3)}{2}$};
		\draw[dashed](\rhoc,\frhoc)--(1.08, \frhoc) node[right] at (1.07,\frhoc) {$ F(\rho_B^3)$};
		\end{tikzpicture}
		
	}
	\caption{Fair merging coupling conditions with $\rho_B^1= 0.2$, $ \rho_B^2= 0.3$ and $\rho_B^3= 0.6$. $\bar \rho$ is the  value found from solving (\ref{eqz}) at the node and $\bar \rho_0$ is the ingoing value  of the solution   of the  layer problems at the node found from the analysis in Section \ref{proof}.}
	\label{couplingtotal}
\end{figure}
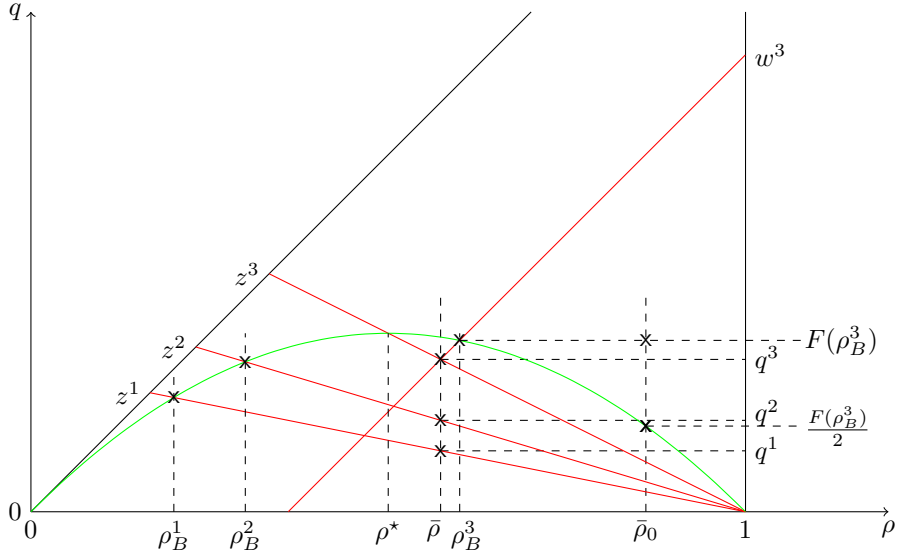

\begin{figure}[h]
	\center
	\externaltikz{solutionnode}{
		\begin{tikzpicture}[scale=0.65]
			\def\rhoa{0.2}
		\def\rhob{0.3}
		\def\rhoc{0.6}
		
		\pgfmathsetmacro{\za}{\rhoa}
		\pgfmathsetmacro{\zb}{\rhob}
		\pgfmathsetmacro{\zc}{\za+\zb}
		\pgfmathsetmacro{\wc}{\rhoc*\rhoc}
		\pgfmathsetmacro{\barrho}{(\wc+\zc)/(1+\zc)}

		\pgfmathsetmacro{\frhoc}{\rhoc*(1-\rhoc)}
	
		\pgfmathsetmacro{\barrhozero}{(1+sqrt(1-2*\frhoc))/2}
	
		\begin{axis}[
		legend style = {at={(1,0)}, xshift=-0.1cm, yshift=0.1cm, anchor=south east},
		legend columns= 2,	
		xlabel = t,		
		ylabel = $\rho$,
		ytick = {\barrho,\barrhozero},
		yticklabels={$\bar \rho$, $\bar \rho_0$},
		]
		\addplot[color = blue!0!red,thick] file{Data/merge_Lindeg_rho_3ex111_eps01_trace.txt};
		\addlegendentry{$\varepsilon=0.1$}
		\addplot[color = blue!33!red,thick] file{Data/merge_Lindeg_rho_3ex111_eps001_trace.txt};
		\addlegendentry{$\varepsilon=0.01$}
		\addplot[color = blue!66!red,thick] file{Data/merge_Lindeg_rho_3ex111_eps0001_trace.txt};
		\addlegendentry{$\varepsilon=0.001$}
		\addplot[color = blue!100!red,thick] file{Data/merge_Lindeg_rho_3ex111_eps00001_trace.txt};		
		\addlegendentry{$\varepsilon=0.0001$}
		\end{axis}
		\end{tikzpicture}
  	}
	\caption{Solution of the relaxation system with  different values of $\varepsilon$ for the initial values in Figure \ref{couplingtotal}. Time development of the density at the junction from $\bar \rho$ to $\bar \rho_0$.}
	\label{junctionvalues}
\centering
\end{figure}
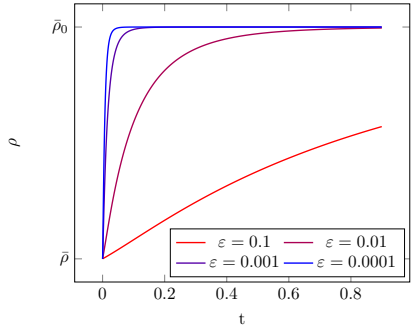

\section{Numerical results}
\label{Numerical results}
In this section we compare  relaxation  and macroscopic network solutions with the  different  coupling conditions for several characteristic numerical examples. 
The relaxation  model is discretized in its conservative form \eqref{eq:lindeg+relax} using a Godunov scheme. A Godunov scheme is also used for  the LWR model. 
In all numerical examples the intervals on the edges $[0,1]$ are discretized with $1000$ cell.
The ingoing edges  are connected to the junction at $x=1$, while the cars enter at $x=0$ into the outgoing edges.
At the outer boundaries zero-Neumann boundary conditions are imposed.
The scaling parameter $\varepsilon $ in the  relaxation system is chosen as $\varepsilon = 0.001$.
As initial conditions the densities $\rho^i$ are chosen constant on each road.
The additional initial condition for $z^i$ in the relaxation  model is chosen in equilibrium $z^i=\frac{F(\rho^i)}{1-\rho^i}$.  
All solutions are computed up to $T=1$.

\subsection{Fair merging}
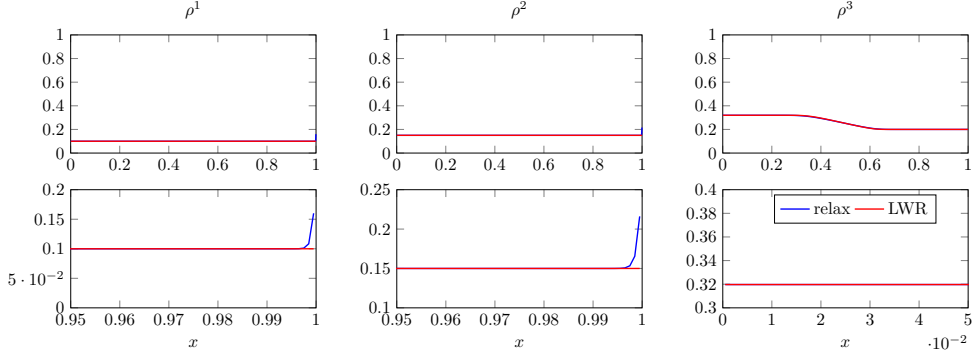
\begin{figure}[h]
	\externaltikz{merge_case11}{
		\begin{tikzpicture}[scale=0.65]
			\begin{groupplot}[
				group style={group size=3 by 2, vertical sep = 0.75cm, horizontal sep = 1.65cm},
				width = 6.6cm,
				height = 4cm,
				xmin = 0.0, xmax = 1.0,
				ymin = 0.0, ymax = 1.0,
				legend style = {at={(0.5,1)},xshift=0.2cm,yshift=-0.1cm,anchor=north},
				legend columns= 3,			
			]
			
			\nextgroupplot[ title = $\rho^1$]
				\addplot[color = blue,thick] file{Data/FairMerge_LindegNew_rho_1ex1_eps0001.txt};
				\addplot[color = red,thick] file{Data/FairMerge_LWR_rho_1ex1.txt};
				
			\nextgroupplot[  title =  $\rho^2$]
				\addplot[color = blue,thick] file{Data/FairMerge_LindegNew_rho_2ex1_eps0001.txt};
				\addplot[color = red,thick] file{Data/FairMerge_LWR_rho_2ex1.txt};
			\nextgroupplot[  title = $\rho^3$]
				\addplot[color = blue,thick] file{Data/FairMerge_LindegNew_rho_3ex1_eps0001.txt};
				\addplot[color = red,thick] file{Data/FairMerge_LWR_rho_3ex1.txt};
			\nextgroupplot[ xlabel = $x$,
				xmin = 0.95, xmax = 1.0,
				ymin = 0.0, ymax = 0.2]	
				\addplot[color = blue,thick] file{Data/FairMerge_LindegNew_rho_1ex1_eps0001.txt};
				\addplot[color = red,thick] file{Data/FairMerge_LWR_rho_1ex1.txt};
			\nextgroupplot[xlabel = $x$,
				xmin = 0.95, xmax = 1.0,
				ymin = 0.1, ymax = 0.25]	
				\addplot[color = blue,thick] file{Data/FairMerge_LindegNew_rho_2ex1_eps0001.txt};
				\addplot[color = red,thick] file{Data/FairMerge_LWR_rho_2ex1.txt};
			\nextgroupplot[ xlabel = $x$,
					xmin = 0.0, xmax = 0.05,
					ymin = 0.3, ymax = 0.4]	
				\addplot[color = blue,thick] file{Data/FairMerge_LindegNew_rho_3ex1_eps0001.txt};
				\addlegendentry{relax}
				\addplot[color = red,thick] file{Data/FairMerge_LWR_rho_3ex1.txt};
				\addlegendentry{LWR}
			\end{groupplot}
		\end{tikzpicture}
	}
	\caption{Fair merging with $\rho^1 = 0.1$, $\rho^2 = 0.15$, $\rho^3 = 0.2$.  Red: solutions of the LWR-equations, blue:  solutions of the relaxation system. First row: solutions on the full domain, second row: zoom around the node.}
	\label{fig:Merge_case1}
\end{figure}

First we compare the numerical solutions of the relaxation model with the coupling conditions  \ref{eqz} with the results obtained for the LWR model with the coupling conditions \eqref{fair}.
In Figure \ref{fig:Merge_case1} the initial densities are chosen as $\rho^1 = 0.1$, $\rho^2 = 0.15$ and  $\rho^3 = 0.2$.
The densities are small enough, such  that all  cars can pass the junction, which corresponds to  Case 7, first subcase. The $\rho_0^i$ are given by $\rho_-(F(\rho_B^1)+F(\rho_B^2))$ with the  numerical value $\rho_0^i=0.3197$.
In Figure \ref{fig:Merge_case1} the numerical solutions are shown. 
Outside of the layer regions, the solution of the relaxation model (blue) is almost identical to the solution of the LWR model (red).
In the second row a zoom into the layer regions at the junctions is shown.
On edge $1$ and $2$ we can observe two boundary layers, as these correspond to stable cases.
In edge $3$ there is no layer, since the half space solution is unstable. The solution at $x=0$ fits exactly to the analytical value
of $\rho_0^i$.

\begin{figure}[h]
	\externaltikz{fairmerge_case2a}{		
		\begin{tikzpicture}[scale=0.65]
			\begin{groupplot}[
				group style={group size=3 by 2, vertical sep = 0.75cm, horizontal sep = 1.75cm},
				width = 6.5cm,
				height = 4cm,
				xmin = -0.0, xmax = 1.0,
				ymin = 0.0, ymax = 1.0,
				legend style = {at={(0.5,1)},xshift=0.2cm,yshift=-0.1cm,anchor=north},
				legend columns= 3,			
				]
			\nextgroupplot[ title = $\rho^1$]
				\addplot[color = blue,thick] file{Data/FairMerge_LindegNew_rho_1ex2_eps0001.txt};
%				\addlegendentry{relax}
				\addplot[color = red,thick] file{Data/FairMerge_LWR_rho_1ex2.txt};
%				\addlegendentry{LWR}
			\nextgroupplot[ title = $\rho^2$]
				\addplot[color = blue,thick] file{Data/FairMerge_LindegNew_rho_2ex2_eps0001.txt};
				\addplot[color = red,thick] file{Data/FairMerge_LWR_rho_2ex2.txt};
			\nextgroupplot[  title = $\rho^3$]
				\addplot[color = blue,thick] file{Data/FairMerge_LindegNew_rho_3ex2_eps0001.txt};
				\addplot[color = red,thick] file{Data/FairMerge_LWR_rho_3ex2.txt};
			\nextgroupplot[ xlabel = $x$,
					xmin = 0.95, xmax = 1.0,
					ymin = 0.7, ymax = 0.9]	
				\addplot[color = blue,thick] file{Data/FairMerge_LindegNew_rho_1ex2_eps0001.txt};
				\addplot[color = red,thick] file{Data/FairMerge_LWR_rho_1ex2.txt};
			\nextgroupplot[ xlabel = $x$,
				xmin = 0.95, xmax = 1.0,
				ymin = 0.7, ymax = 0.9]	
				\addplot[color = blue,thick] file{Data/FairMerge_LindegNew_rho_2ex2_eps0001.txt};
				\addplot[color = red,thick] file{Data/FairMerge_LWR_rho_2ex2.txt};
			\nextgroupplot[ xlabel = $x$,
				xmin = 0.0, xmax = 0.05,
				ymin = 0.4, ymax = 0.9]	
				\addplot[color = blue,thick] file{Data/FairMerge_LindegNew_rho_3ex2_eps0001.txt};
				\addlegendentry{relax}
				\addplot[color = red,thick] file{Data/FairMerge_LWR_rho_3ex2.txt};
				\addlegendentry{LWR}
			\end{groupplot}
		\end{tikzpicture}
	}
	\caption{Fair merging with $\rho^1 = 0.7$, $\rho^2 = 0.6$, $\rho^3 = 0.2$.  }
	\label{fig:Merge_case2}
\end{figure}
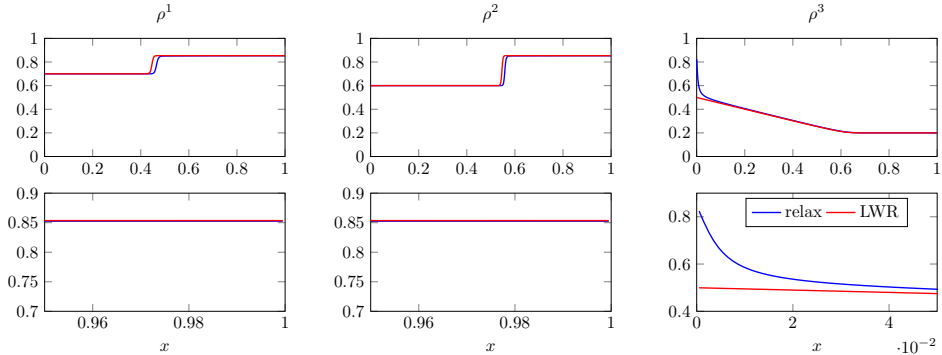

In Figure \ref{fig:Merge_case2} the numerical solutions to the initial values $\rho^1 = 0.7$, $\rho^2 = 0.6$ and  $\rho^3 = 0.2$ are shown. 
In this situation more cars are approaching the junction than can enter  road $3$. We are in the situation of Case 1 with the analytical value $\rho_0^i =\rho_+(\sigma/2)=0.83536$.
Thus the flow in the exiting road is set to its maximum, while there are jams propagating upstream in the ingoing roads.
Here we observe only in edge $3$ a layer, which interacts with the tail of the rarefaction wave.
In the ingoing roads the unstable layer solution enforce the new values at the junction.
In these roads the shock waves of the relaxation model are slightly behind those of the macroscopic one. 
This stems from an initial layer, as the layer at the junction has to form at the beginning, see Figure \ref{junctionvalues}. 
This happens in short time and is not visible at the rarefaction waves, but it remains noticeable at the shocks.
The speeds of the shocks is identical in both models, as the connected states coincide, i.e. the delay does not change over time.

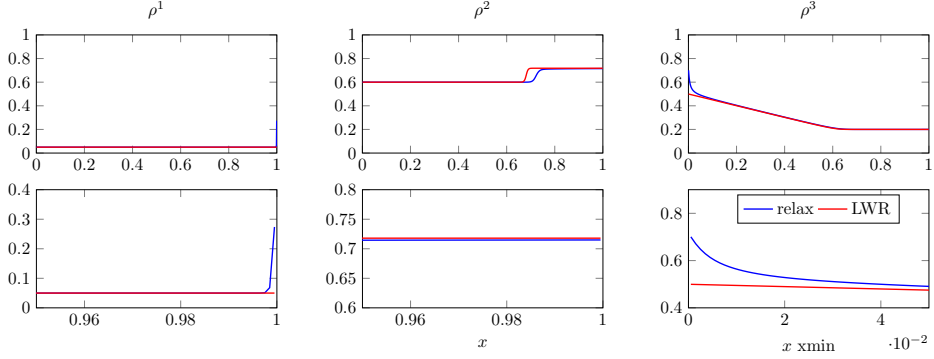
\begin{figure}[h]
	\externaltikz{fairmerge_case3}{
		\begin{tikzpicture}[scale=0.65]
			\begin{groupplot}[
				group style={group size=3 by 2, vertical sep = 0.75cm, horizontal sep = 1.75cm},
				width = 6.5cm,
				height = 4cm,
				xmin = -0.0, xmax = 1.0,
				ymin = 0.0, ymax = 1.0,
				legend style = {at={(0.5,1)},xshift=0.2cm,yshift=-0.1cm,anchor=north},
				legend columns= 3,			
				]		
			\nextgroupplot[ title = $\rho^1$]
				\addplot[color = blue,thick] file{Data/FairMerge_LindegNew_rho_1ex3_eps0001.txt};
%				\addlegendentry{relax}
				\addplot[color = red,thick] file{Data/FairMerge_LWR_rho_1ex3.txt};
%				\addlegendentry{LWR}
			\nextgroupplot[ title = $\rho^2$]
				\addplot[color = blue,thick] file{Data/FairMerge_LindegNew_rho_2ex3_eps0001.txt};
				\addplot[color = red,thick] file{Data/FairMerge_LWR_rho_2ex3.txt};
			\nextgroupplot[  title = $\rho^3$]
				\addplot[color = blue,thick] file{Data/FairMerge_LindegNew_rho_3ex3_eps0001.txt};
				\addplot[color = red,thick] file{Data/FairMerge_LWR_rho_3ex3.txt};
			\nextgroupplot[ %ylabel = $\rho^1$,
				xmin = 0.95, xmax = 1.0,
				ymin = 0.0, ymax = 0.4]	
				\addplot[color = blue,thick] file{Data/FairMerge_LindegNew_rho_1ex3_eps0001.txt};
				\addplot[color = red,thick] file{Data/FairMerge_LWR_rho_1ex3.txt};
			\nextgroupplot[ xlabel = $x$,
				xmin = 0.95, xmax = 1.0,
				ymin = 0.6, ymax = 0.8]	
				\addplot[color = blue,thick] file{Data/FairMerge_LindegNew_rho_2ex3_eps0001.txt};
				\addplot[color = red,thick] file{Data/FairMerge_LWR_rho_2ex3.txt};
			\nextgroupplot[ xlabel = $x$
				xmin = 0.0, xmax = 0.05,
				ymin = 0.4, ymax = 0.9]	
				\addplot[color = blue,thick] file{Data/FairMerge_LindegNew_rho_3ex3_eps0001.txt};
				\addlegendentry{relax}
				\addplot[color = red,thick] file{Data/FairMerge_LWR_rho_3ex3.txt};
				\addlegendentry{LWR}
			\end{groupplot}
		\end{tikzpicture}
	}
	\caption{Fair merging with $\rho^1 = 0.05$, $\rho^2 = 0.6$, $\rho^3 = 0.2$. }
	\label{fig:Merge_case3}
\end{figure}

In the next example, with the initial values $\rho^1 = 0.05$, $\rho^2 = 0.6$ and $\rho^3 = 0.2$, few cars enter from road $1$ but many from road $2$. We are in Case 4, first subcase. The analytical value at the junction is $\rho_0^i=\rho_+(\sigma-F(\rho_B^1))=0.7179$.
As shown in Figure \ref{fig:Merge_case3}, the flow in road $3$ is at maximum such that all cars from road $1$ and most of road $2$ can pass.
Those which do not fit in, create a jam in road $2$.
Again we see a delay of the shock, as in the previous example. 
Similarly we observe a layer in edge $3$.
But here also a layer in road $1$ is present, as the solution of the half space is now stable.

\begin{remark}
We mention that the numerical investigation of conditions (\ref{mixed2})	gives slightly different  values for the relaxation system at the nodes, but the same results in the interior of the domain. That means,  also in this case, the relaxation system leads  to the LWR equations  on the network with the fair merging condition (\ref{fair}).
\end{remark}
\subsection{Merging with priority lane}
Here,  the numerical solutions of the relaxation  model with the coupling conditions \ref{eq:merge_4} are compared to those obtained for the LWR model with the coupling conditions  \ref{prio}.

In the first example with $\rho^1 = 0.6$, $\rho^2 = 0.7$ and $\rho^3 = 0.2$, shown in Figure \ref{fig:Merge_case10}, many cars arrive at the junction. 
\begin{figure}[h]
	\externaltikz{merge_case10}{
		\begin{tikzpicture}[scale=0.65]
		\begin{groupplot}[
		group style={group size=3 by 2, vertical sep = 0.75cm, horizontal sep = 1.75cm},
		width = 6.5cm,
		height = 4cm,
		xmin = -0.0, xmax = 1.0,
		ymin = 0.0, ymax = 1.0,
		legend style = {at={(0.5,1)},xshift=0.2cm,yshift=-0.1cm,anchor=north},
		legend columns= 3,			
		]		
		\nextgroupplot[ title = $\rho^1$]
		\addplot[color = blue,thick] file{Data/PriorityMerge_LindegNew_rho_1ex10_eps0001.txt};
%		\addlegendentry{relax}
		\addplot[color = red,thick] file{Data/PriorityMerge_LWR_rho_1ex10.txt};
%		\addlegendentry{LWR}
		\nextgroupplot[ title = $\rho^2$]
		\addplot[color = blue,thick] file{Data/PriorityMerge_LindegNew_rho_2ex10_eps0001.txt};
		\addplot[color = red,thick] file{Data/PriorityMerge_LWR_rho_2ex10.txt};
		\nextgroupplot[  title = $\rho^3$]
		\addplot[color = blue,thick] file{Data/PriorityMerge_LindegNew_rho_3ex10_eps0001.txt};
		\addplot[color = red,thick] file{Data/PriorityMerge_LWR_rho_3ex10.txt};
		\nextgroupplot[ xlabel = $x$, %ylabel = $\rho^1$,
		xmin = 0.95, xmax = 1.0,
		ymin = 0.2, ymax = 0.6]	
		\addplot[color = blue,thick] file{Data/PriorityMerge_LindegNew_rho_1ex10_eps0001.txt};
		\addplot[color = red,thick] file{Data/PriorityMerge_LWR_rho_1ex10.txt};
		\nextgroupplot[ xlabel = $x$,  %ylabel = $\rho^2$,
		xmin = 0.95, xmax = 1.0,
		ymin = 0.9, ymax = 1.1]	
		\addplot[color = blue,thick] file{Data/PriorityMerge_LindegNew_rho_2ex10_eps0001.txt};
		\addplot[color = red,thick] file{Data/PriorityMerge_LWR_rho_2ex10.txt};
		\nextgroupplot[ xlabel = $x$, %  ylabel = $\rho^3$]
		xmin = 0.0, xmax = 0.05,
		ymin = 0.45, ymax = 0.8]	
		\addplot[color = blue,thick] file{Data/PriorityMerge_LindegNew_rho_3ex10_eps0001.txt};
		\addlegendentry{relax}
		\addplot[color = red,thick] file{Data/PriorityMerge_LWR_rho_3ex10.txt};
		\addlegendentry{LWR}
		\end{groupplot}
		\end{tikzpicture}
	}
	\caption{Priority merge with $\rho^1 = 0.6$, $\rho^2 = 0.7$, $\rho^3 = 0.2$. }
	\label{fig:Merge_case10}
\end{figure}
As those of road $1$ have priority, the maximal flow is established, while all cars in road $2$ have to wait. 
Layers can be observed in road $1$ and $3$.

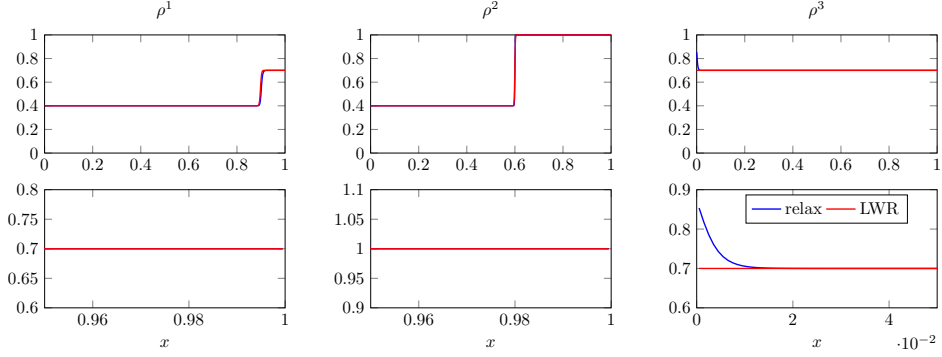
\begin{figure}[h]
	\externaltikz{merge_case12}{
		\begin{tikzpicture}[scale=0.65]
		\begin{groupplot}[
		group style={group size=3 by 2, vertical sep = 0.75cm, horizontal sep = 1.75cm},
		width = 6.5cm,
		height = 4cm,
		xmin = -0.0, xmax = 1.0,
		ymin = 0.0, ymax = 1.0,
		legend style = {at={(0.5,1)},xshift=0.2cm,yshift=-0.1cm,anchor=north},
		legend columns= 3,			
		]		
		\nextgroupplot[ title = $\rho^1$]
		\addplot[color = blue,thick] file{Data/PriorityMerge_LindegNew_rho_1ex12_eps0001.txt};
%		\addlegendentry{relax}
		\addplot[color = red,thick] file{Data/PriorityMerge_LWR_rho_1ex12.txt};
%		\addlegendentry{LWR}
		\nextgroupplot[ title = $\rho^2$]
		\addplot[color = blue,thick] file{Data/PriorityMerge_LindegNew_rho_2ex12_eps0001.txt};
		\addplot[color = red,thick] file{Data/PriorityMerge_LWR_rho_2ex12.txt};
		\nextgroupplot[  title = $\rho^3$]
		\addplot[color = blue,thick] file{Data/PriorityMerge_LindegNew_rho_3ex12_eps0001.txt};
		\addplot[color = red,thick] file{Data/PriorityMerge_LWR_rho_3ex12.txt};
		\nextgroupplot[ xlabel = $x$, %ylabel = $\rho^1$,
		xmin = 0.95, xmax = 1.0,
		ymin = 0.6, ymax = 0.8]	
		\addplot[color = blue,thick] file{Data/PriorityMerge_LindegNew_rho_1ex12_eps0001.txt};
		\addplot[color = red,thick] file{Data/PriorityMerge_LWR_rho_1ex12.txt};
		\nextgroupplot[ xlabel = $x$, %ylabel = $\rho^2$,
		xmin = 0.95, xmax = 1.0,
		ymin = 0.9, ymax = 1.1]	
		\addplot[color = blue,thick] file{Data/PriorityMerge_LindegNew_rho_2ex12_eps0001.txt};
		\addplot[color = red,thick] file{Data/PriorityMerge_LWR_rho_2ex12.txt};
		\nextgroupplot[ xlabel = $x$,%  ylabel = $\rho^3$]
		xmin = 0.0, xmax = 0.05,
		ymin = 0.6, ymax = 0.9]	
		\addplot[color = blue,thick] file{Data/PriorityMerge_LindegNew_rho_3ex12_eps0001.txt};
		\addlegendentry{relax}
		\addplot[color = red,thick] file{Data/PriorityMerge_LWR_rho_3ex12.txt};
		\addlegendentry{LWR}
		\end{groupplot}
		\end{tikzpicture}
	}
	\caption{Priority merge with $\rho^1 = 0.4$, $\rho^2 = 0.4$, $\rho^3 = 0.7$. }
	\label{fig:Merge_case12}
\end{figure}

In the second example we consider a situation, with  the same amount of cars  in the ingoing roads and only little space in the outgoing one
$\rho^1 = 0.4$, $\rho^2 = 0.4$, $\rho^3 = 0.7$.
As expected, we can see in Figure \ref{fig:Merge_case12} that all the cars in road $2$ have to wait and thus a larger shock forms.
Not all cars in the first road can pass, but the flow is larger as in the second road. 
%This is again a situation from Case B.

In both cases we obtain a perfect numerical agreement of relaxation solutions and LWR-solutions on the network
outside of the layers.

\subsection{Diverging with drivers preferences}

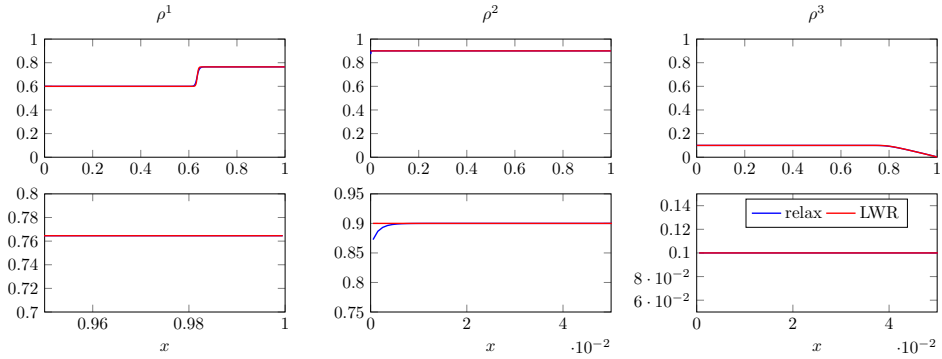
\begin{figure}[h!]
	\externaltikz{split_case8}{
		\begin{tikzpicture}[scale=0.65]
		\begin{groupplot}[
		group style={group size=3 by 2, vertical sep = 0.75cm, horizontal sep = 1.75cm},
		width = 6.5cm,
		height = 4cm,
		xmin = -0.0, xmax = 1.0,
		ymin = 0.0, ymax = 1.0,
		legend style = {at={(0.5,1)},xshift=0.2cm,yshift=-0.1cm,anchor=north},
		legend columns= 3,			
		]
		
		\nextgroupplot[ title = $\rho^1$]
		\addplot[color = blue,thick] file{Data/DivergeAlpha_LindegNew_rho_1ex21_eps0001.txt};
		\addplot[color = red,thick] file{Data/DivergeAlpha_LWR_rho_1ex21.txt};
		
		\nextgroupplot[ title = $\rho^2$]
		\addplot[color = blue,thick] file{Data/DivergeAlpha_LindegNew_rho_2ex21_eps0001.txt};
		\addplot[color = red,thick] file{Data/DivergeAlpha_LWR_rho_2ex21.txt};
		
		\nextgroupplot[ title = $\rho^3$]
		\addplot[color = blue,thick] file{Data/DivergeAlpha_LindegNew_rho_3ex21_eps0001.txt};
		\addplot[color = red,thick] file{Data/DivergeAlpha_LWR_rho_3ex21.txt};
		\nextgroupplot[ xlabel =  $x$, 
		xmin = 0.95, xmax = 1.0,
		ymin = 0.7, ymax = 0.8]	
		\addplot[color = blue,thick] file{Data/DivergeAlpha_LindegNew_rho_1ex21_eps0001.txt};
		\addplot[color = red,thick] file{Data/DivergeAlpha_LWR_rho_1ex21.txt};
		\nextgroupplot[ xlabel =  $x$, 
		xmin = 0.0, xmax = 0.05,
		ymin = 0.75, ymax = 0.95]	
		\addplot[color = blue,thick] file{Data/DivergeAlpha_LindegNew_rho_2ex21_eps0001.txt};
		\addplot[color = red,thick] file{Data/DivergeAlpha_LWR_rho_2ex21.txt};
		\nextgroupplot[ xlabel =  $x$, 
		xmin = 0.0, xmax = 0.05,
		ymin = 0.05, ymax = 0.15]	
		\addplot[color = blue,thick] file{Data/DivergeAlpha_LindegNew_rho_3ex21_eps0001.txt};
		\addlegendentry{relax}
		\addplot[color = red,thick] file{Data/DivergeAlpha_LWR_rho_3ex21.txt};
		\addlegendentry{LWR}
		\end{groupplot}
		
		\end{tikzpicture}
	}
	\caption{Diverging with driver preferences: $\rho^1 = 0.6$, $\rho^2 = 0.9$, $\rho^3 = 0.0$ }
	\label{fig:Split_case8}
\end{figure}

We consider a  junction with equal drivers preferences $\alpha=\frac{1}{2}$ and conditions \ref{divmacro1}  for the macroscopic equations and (\ref{eq:div_2}), (\ref{pref}) for the relaxation problem. We consider two   examples. 	
First, with the initial conditions $\rho^1 = 0.8$, $\rho^2 = 0.1$ and $\rho^3 = 0.3$, there is enough space in both outgoing roads such that the maximal flow can be established, as shown in Figure \ref{fig:Split_case7}. 
%We have $C^1=c^1$.
\begin{figure}[h!]
	\externaltikz{split_case7}{
		\begin{tikzpicture}[scale=0.65]
		\begin{groupplot}[
		group style={group size=3 by 2, vertical sep = 0.75cm, horizontal sep = 1.75cm},
		width = 6.5cm,
		height = 4cm,
		xmin = -0.0, xmax = 1.0,
		ymin = 0.0, ymax = 1.0,
		legend style = {at={(0.5,1)},xshift=0.2cm,yshift=-0.1cm,anchor=north},
		legend columns= 3,			
		]
		
		\nextgroupplot[ title = $\rho^1$]
		\addplot[color = blue,thick] file{Data/DivergeAlpha_LindegNew_rho_1ex20_eps0001.txt};
%		\addlegendentry{relax}
		\addplot[color = red,thick] file{Data/DivergeAlpha_LWR_rho_1ex20.txt};
		\nextgroupplot[ title = $\rho^2$]
		\addplot[color = blue,thick] file{Data/DivergeAlpha_LindegNew_rho_2ex20_eps0001.txt};
		\addplot[color = red,thick] file{Data/DivergeAlpha_LWR_rho_2ex20.txt};
		\nextgroupplot[ title = $\rho^3$]
		\addplot[color = blue,thick] file{Data/DivergeAlpha_LindegNew_rho_3ex20_eps0001.txt};
		\addplot[color = red,thick] file{Data/DivergeAlpha_LWR_rho_3ex20.txt};
		\nextgroupplot[ xlabel =  $x$, 
		xmin = 0.95, xmax = 1.0,
		ymin = 0.25, ymax = 0.55]	
		\addplot[color = blue,thick] file{Data/DivergeAlpha_LindegNew_rho_1ex20_eps0001.txt};
		\addplot[color = red,thick] file{Data/DivergeAlpha_LWR_rho_1ex20.txt};
		\nextgroupplot[ xlabel =  $x$, 
		xmin = 0.0, xmax = 0.05,
		ymin = 0.1, ymax = 0.2]	
		\addplot[color = blue,thick] file{Data/DivergeAlpha_LindegNew_rho_2ex20_eps0001.txt};
		\addplot[color = red,thick] file{Data/DivergeAlpha_LWR_rho_2ex20.txt};
		\nextgroupplot[ xlabel =  $x$, 
		xmin = 0.0, xmax = 0.05,
		ymin = 0.1, ymax = 0.2]	
		\addplot[color = blue,thick] file{Data/DivergeAlpha_LindegNew_rho_3ex20_eps0001.txt};
		\addlegendentry{relax}
		\addplot[color = red,thick] file{Data/DivergeAlpha_LWR_rho_3ex20.txt};
		\addlegendentry{LWR}
		\end{groupplot}
		
		\end{tikzpicture}
	}
	\caption{Diverging with driver preferences: $\rho^1 = 0.8$, $\rho^2 = 0.1$, $\rho^3 = 0.3$ }
	\label{fig:Split_case7}
\end{figure}
In Figure \ref{fig:Split_case8} the solutions corresponding to the initial values $\rho^1 = 0.6$, $\rho^2 = 0.9$ and $\rho^3 = 0.0$ are shown.
Although road $3$ is completely free, only few case can enter, as their way is blocked by cars waiting to enter road $2$. 
%In this case $C^1=2 c^2$.
Thus the high density on road $2$ is causing a left going shock on the ingoing road.
A layer forms only on road $2$, since on the other two the macroscopic characteristics move away from the junction.

\subsection{Diverging without preferences}

\begin{figure}[h!]
	\externaltikz{split_case4}{
		\begin{tikzpicture}[scale=0.65]
		\begin{groupplot}[
				group style={group size=3 by 2, vertical sep = 0.75cm, horizontal sep = 1.75cm},
				width = 6.5cm,
				height = 4cm,
				xmin = -0.0, xmax = 1.0,
				ymin = 0.0, ymax = 1.0,
				legend style = {at={(0.5,1)},xshift=0.2cm,yshift=-0.1cm,anchor=north},
				legend columns= 3,			
				]
		\nextgroupplot[ title = $\rho^1$]
			\addplot[color = blue,thick] file{Data/FreeDiverge_LindegNew_rho_1ex40_eps0001.txt};
%			\addlegendentry{relax}
			\addplot[color = red,thick] file{Data/FreeDiverge_LWR_rho_1ex40.txt};
%			\addlegendentry{LWR}
		\nextgroupplot[ title = $\rho^2$]
			\addplot[color = blue,thick] file{Data/FreeDiverge_LindegNew_rho_2ex40_eps0001.txt};
			\addplot[color = red,thick] file{Data/FreeDiverge_LWR_rho_2ex40.txt};
		
		\nextgroupplot[ title = $\rho^3$]
			\addplot[color = blue,thick] file{Data/FreeDiverge_LindegNew_rho_3ex40_eps0001.txt};
			\addplot[color = red,thick] file{Data/FreeDiverge_LWR_rho_3ex40.txt};
		\nextgroupplot[
				xmin = 0.95, xmax = 1.0,
				ymin = 0.25, ymax = 0.6]	
			\addplot[color = blue,thick] file{Data/FreeDiverge_LindegNew_rho_1ex40_eps0001.txt};
			\addplot[color = red,thick] file{Data/FreeDiverge_LWR_rho_1ex40.txt};
		\nextgroupplot[ xlabel =  $x$, 
				xmin = 0.0, xmax = 0.05,
				ymin = 0.1, ymax = 0.2]	
			\addplot[color = blue,thick] file{Data/FreeDiverge_LindegNew_rho_2ex40_eps0001.txt};
			\addplot[color = red,thick] file{Data/FreeDiverge_LWR_rho_2ex40.txt};
		\nextgroupplot[ xlabel =  $x$, 
				xmin = 0.0, xmax = 0.05,
				ymin = 0.1, ymax = 0.2]	
			\addplot[color = blue,thick] file{Data/FreeDiverge_LindegNew_rho_3ex40_eps0001.txt};
			\addlegendentry{relax}
			\addplot[color = red,thick] file{Data/FreeDiverge_LWR_rho_3ex40.txt};
			\addlegendentry{LWR}
		\end{groupplot}
		
		\end{tikzpicture}
	}
	\caption{Diverging without preferences: $\rho^1 = 0.7$, $\rho^2 = 0.2$, $\rho^3 = 0.1$. }
	\label{fig:Split_case4}
\end{figure}
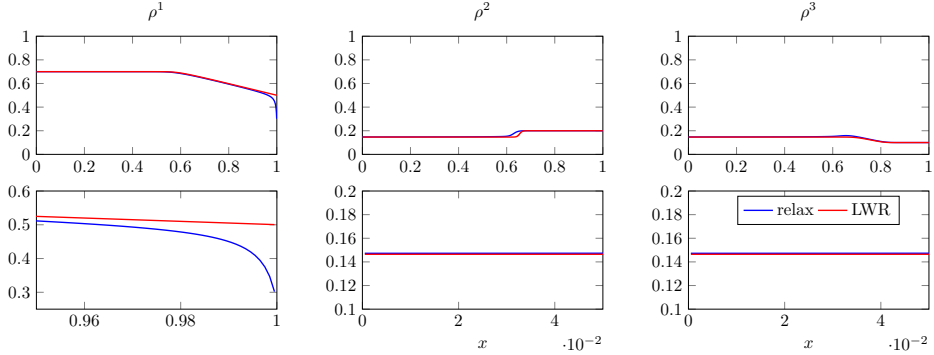

We consider the situation with condition \ref{divmacro2} for the macroscopic and conditions  \ref{eq:div_3}, \ref{ex2}
for the relaxation model.
The first example investigates the case $\rho^1 = 0.7$, $\rho^2 = 0.2$, $\rho^3 = 0.1$, see Figure
\ref{fig:Split_case4}.
In this case, the  flux is distributed equally onto the outgoing roads, such that a small shock and a small rarefaction wave arise. 
On the right hand side we observe that a  layer forms in the first road but not in the two exiting ones.   

\begin{figure}[h!]
	\externaltikz{split_case6}{
		\begin{tikzpicture}[scale=0.65]
		\begin{groupplot}[
		group style={group size=3 by 2, vertical sep = 0.75cm, horizontal sep = 1.75cm},
		width = 6.5cm,
		height = 4cm,
		xmin = -0.0, xmax = 1.0,
		ymin = 0.0, ymax = 1.0,
		legend style = {at={(0.5,0)},xshift=0.2cm,yshift=0.1cm,anchor=south},
		legend columns= 3,			
		]
		
		\nextgroupplot[ title = $\rho^1$]
		\addplot[color = blue,thick] file{Data/FreeDiverge_LindegNew_rho_1ex42_eps0001.txt};
		
%		\addlegendentry{kinetic}
		\addplot[color = red,thick] file{Data/FreeDiverge_LWR_rho_1ex42.txt};
%		\addlegendentry{LWR}
		
		\nextgroupplot[ title = $\rho^2$]
		\addplot[color = blue,thick] file{Data/FreeDiverge_LindegNew_rho_2ex42_eps0001.txt};
		\addplot[color = red,thick] file{Data/FreeDiverge_LWR_rho_2ex42.txt};
		
		\nextgroupplot[ title = $\rho^3$]
		\addplot[color = blue,thick] file{Data/FreeDiverge_LindegNew_rho_3ex42_eps0001.txt};
		\addplot[color = red,thick] file{Data/FreeDiverge_LWR_rho_3ex42.txt};
		
		\nextgroupplot[ xlabel =  $x$,
		xmin = 0.95, xmax = 1.0,
		ymin = 0.25, ymax = 0.55]	
		\addplot[color = blue,thick] file{Data/FreeDiverge_LindegNew_rho_1ex42_eps0001.txt};
		\addplot[color = red,thick] file{Data/FreeDiverge_LWR_rho_1ex42.txt};
		
		\nextgroupplot[ xlabel =  $x$,
		xmin = 0.0, xmax = 0.05,
		ymin = 0.25, ymax = 0.35]	
		\addplot[color = blue,thick] file{Data/FreeDiverge_LindegNew_rho_2ex42_eps0001.txt};
		\addplot[color = red,thick] file{Data/FreeDiverge_LWR_rho_2ex42.txt};

		\nextgroupplot[ xlabel =  $x$, 
		xmin = 0.0, xmax = 0.05,
		ymin = 0.3, ymax = 1]	
		\addplot[color = blue,thick] file{Data/FreeDiverge_LindegNew_rho_3ex42_eps0001.txt};
		\addlegendentry{kinetic}
		\addplot[color = red,thick] file{Data/FreeDiverge_LWR_rho_3ex42.txt};
		\addlegendentry{LWR}
		\end{groupplot}
		
		\end{tikzpicture}
	}
	\caption{Diverging without preferences: $\rho^1 = 0.6$, $\rho^2 = 0.1$, $\rho^3 = 0.95$ }
	\label{fig:Split_case6}
\end{figure}
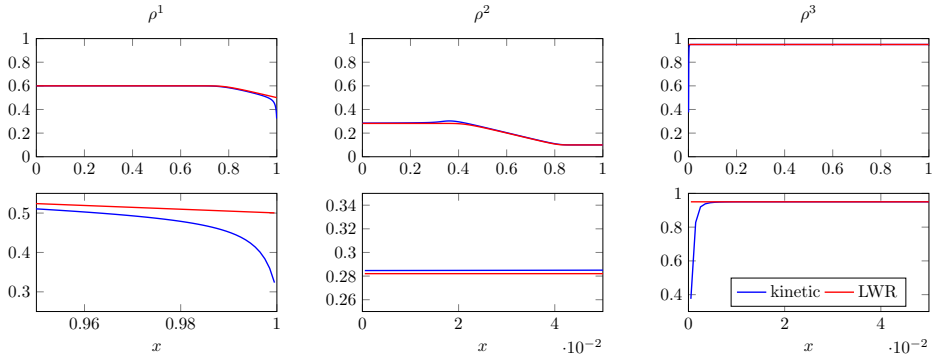

In the second example, Figure \ref{fig:Split_case6}, the results with the initial conditions $\rho^1 = 0.6$, $\rho^2 = 0.1$ and $\rho^3 = 0.95$ are shown.
Since the traffic on road $3$ is dense only very few cars enter there. 
Most of the vehicles  enter into  road $2$.
In the solution of the relaxation model we  observe two layers, one interacting with the rarefaction wave on road $1$ and one due to the ingoing characteristics on road $3$.
As for the merging case, we observe a very good numerical agreement between the solutions of the relaxation
model and the LWR solution away from the layers at the nodes.

\section{Conclusions}
	We have introduced  general  coupling conditions for a TLD-relaxation model for LWR-networks
	These coupling conditions are related, via  an asymptotic analysis  at the nodes   to well known coupling conditions for the LWR-network.
	
	The  asymptotic analysis  shows that a classical merge condition  for the  nonlinear scalar conservation law
	is related  in the zero-relaxation limit to  an equal density  coupling condition for the relaxation system.
	The numerical findings support and illustrate  the analytical results. One also observes  that there is a range of coupling conditions
	for the relaxation model leading to the same merge  condition  for the scalar conservation law.
	
	For the case of a diverging junction coupling conditions have been defined respecting the physical invariant domain.
	They are investigated numerically  on the network in the limit as $\epsilon \rightarrow 0 $ showing  
 again agreement of the 
	numerical solutions of the relaxation system and the LWR solution for small values of $\epsilon $ .
	
	Finally, we remark, that  the analytical procedure, presented here for the merging case and a special coupling condition, could be extended   to the diverging  case or other coupling conditions in the merging case.

%------------------------------------------------------------------------------
\end{document}